\documentstyle[12pt]{article}
\def\Bbb R{{\rm \bf R}}
\def\proclaim#1{\vskip2mm{\bf #1}\em}
\def\endproclaim{\em \vskip2mm}
\def\tag#1{\eqno(#1)}
\def\gathered{\begin{array}{c}}
\def\endgathered{\end{array}}
\def\text{\mbox}

\begin{document}

\title {Integrating the probe and singular sources methods}
\author{Masaru IKEHATA
\footnote{
Laboratory of Mathematics,
Graduate School of Advanced Science and Engineering,
Hiroshima University, Higashihiroshima 739-8527, JAPAN.
e-mail address: ikehataprobe@gmail.com}
\footnote{Professor Emeritus at Gunma University}
\footnote{Professor Emeritus at Hiroshima University}
}
\maketitle
\begin{abstract}
The probe and singular sources methods are two well-known classical direct reconstruction methods
in inverse obstacle problems governed by partial differential equations.
In this paper, by considering an inverse obstacle problem governed by the Laplace equation in a bounded domain as a prototype case,
an integrated theory of the probe and singular sources methods is proposed.
The theory consists of  three parts: (i) introducing the singular sources method combined with the notion of 
the probe method;
(ii)  finding {\it a third indicator function} whose two ways decomposition yields the indicator functions in the probe and singular sources methods;
(iii) finding the completely integrated version of the probe and singular sources methods.


\noindent
AMS: 35R30, 35J05, 78A46


\end{abstract}


\section{Introduction}

In 1996-2000, 
the area of the mathematical study of the so-called {\it inverse obstacle problems} 
governed by partial differential equations drastically changed
the direction of the research.  It was the appearance of  a group of so-called direct methods, or analytical methods
which is different from that of the traditional optimization methods together with linearization procedures. 
The group consists of five classical methods, that is, the linear sampling method \cite{CKi} of Colton-Kirsch (1996),
the probe method  \cite{IProbe,I2, Iwave} of Ikehata (1998, 1999), the factorization method \cite{K, K2, K3} of Kirsch (1998,1999, 2000), the enclosure method \cite{I1, E00} of Ikehata (1999, 2000) and the singular sources method \cite{P1} of Potthast (2000).
Since then, the range of applications of those classical methods has expanded dramatically.
Today, there is no question that the development of those methods has inspired the search for other direct methods,
such as the no-response test \cite{LU}, the monotonicity method \cite{TR,HU}, etc.

This paper focuses on the relationship between the probe and singular sources methods.
Ignoring a difference in their original presentations, one sees that
the both methods have the common style.  They calculate
a ``field'' defined outside the unknown obstacle from the infinitely many observation data.
The field was called the indicator function firstly in the probe method, which blows up 
on the surface of the obstacle.  
From this property one can extract information about the geometry of the obstacle.

The original indicator function for the singular sources method, defined in the exterior domain of an obstacle,
is given by the scattered field of the spherical wave generated by a point source placed outside the obstacle.
More precisely, it is given by the value of the scattered wave at the location of the point source
and the singular sources method is based on its blowing up property on the surface of the obstacle, see Lemmata 12 and 16 in  \cite{P1}.
The choice was reasonable since the uniqueness proof of  Kirsch-Kress \cite{KiKr} employed the blowing up property
of the field in the case of a sound-soft obstacle or its Neumann derivative in the case of a sound-hard obstacle
on the surface of the obstacle when the point source approaches the point.
In that case, it may be a natural choice.

The indicator function for the probe method first appeared a decade after the publication of 
Isakov's uniqueness theorem \cite{Is} on a special version of the Calder\'on problem \cite{Cal}.
It is given by taking a limit of the difference of two Dirichlet-to-Neumann maps acting on sequences of specially chosen Dirichlet data\footnote{Or Neumann-to-Dirichlet maps acting on common sequences of specially chosen Neumann data, depending on applications.}.
Those two maps are associated with 
the governing equations of the used signal, one is coming from the given medium to be probed and another is coming 
from the background medium possible to be calculated in advance.  This indicator function 
is also defined outside an unknown inclusion in \cite{IProbe, IProbe2} or obstacle in \cite{Iwave} and has the same blowing up property as that of 
the singular sources method.

The problems firstly discussed in both methods are different from each other, however, they have the common idea mentioned above.
This suggests us some kind of a relationship between two methods.
For this,  Nakamura-Potthast-Sini in \cite{NPS} and Nakamura-Potthast in 
\cite{NP} have published an article and a book, respectively, in which they considered
the singular sources method applied to a bounded domain case and compared the probe method
with this singular sources method.
More precisely, they introduced an indicator function for the singular sources method in a bounded domain and
stated that the indicator functions in both methods
coincides with modulo $O(1)$ as the common independent variables of the both functions
approaches a point on the surface of unknown obstacle.
This implies that the blowing up profile of the two indicator functions coincides with each other.
By this,  they concluded that both methods are equivalent to each other despite of the difference of the calculation methods for the indicator functions.  
In addition, they stated also that the justification of the  no-response test (of multi waves as they call)  can be derived
from both the probe and singular sources methods.  However, it should be pointed out that their formulation
of the probe method is given {\it in terms of the singular sources method}.

In this paper, we further pursuit a relationship between the probe and singular sources methods
from a different point of view.  That is the creation of an integrated theory of the probe and singular sources methods.
For the purpose, in contrast to Nakamura-Potthast\cite{NP} and Nakamura-Potthast-Sini \cite{NPS} we formulate the singular sources method
in a bounded domain {\it in terms of the notion of the probe method}.  This makes the description about the singular sources method
so simple and its comparison with the probe method easier.
Then, we found that there is a {\it third indicator function} which blows up on both the outer surface of the given domain  and the surface of an obstacle embedded therein.
The indicator functions  in the probe and singular source methods
are derived as a result of decomposing  the third indicator function  into two ways.  
Besides, by this decomposition of the third indicator function, one can clearly understand
the blowing up property of the indicator function of the singular sources method.
This aspect has not been known before and by this, the author thinks, finally, 
the true meaning of the indicator functions for both methods was understood.
This is an unexpected discovery beyond the previous relationship.
And as an unexpected byproduct, we found the completely integrated version of both methods.  
The method has not only the Side A but also Side B which is previously known only for the probe method \cite{IHokkai}.

\subsection{Setting up the prototype problem}

To discuss the essential part we consider only the simplest and typical mathematical model described below.
The author thinks that it is a source of various reconstruction methods and also the best problem to develop 
or try different methods.

Let $\Omega$ be a bounded domain of $\Bbb R^3$ with Lipschitz boundary.
Assume that $D$ is an open subset of $\Bbb R^3$ with Lipschitz boundary, $\overline{D}\subset\Omega$ and that $\Omega\setminus\overline{D}$ is connected.  This is the starting assumption.

Given $f\in H^{\frac{1}{2}}(\partial\Omega)$ let $u\in H^1(\Omega\setminus\overline{D})$ be the weak solution of
$$
\left\{
\begin{array}{ll}
\displaystyle
\Delta u=0, & x\in\Omega\setminus\overline{D},
\\
\\
\displaystyle
\frac{\partial u}{\partial\nu}=0, & x\in\partial D,
\\
\\
\displaystyle
u=f, & x\in\partial\Omega.
\end{array}
\right.
\tag {1.1}
$$
Define 
$$\begin{array}{ll}
\displaystyle
\Lambda_Df=\frac{\partial u}{\partial\nu}\vert_{\partial\Omega}, & x\in\partial\Omega.
\end{array}
$$
We denote $\Lambda_D$ by $\Lambda_0$ if $D=\emptyset$.

We consider the probe and singular sources method for the following problem.

$\quad$

{\noindent\bf Problem.}  Extract information about the geometry of $D$ from the Dirichlet-to-Neumann map
$\Lambda_D$.

$\quad$

Note that, in Nakamura-Potthast \cite{NP},  the homogeneous Dirichlet boundary condition case $u=0$ on $\partial D$ 
under the assumption that the governing equation is given by the Helmholtz equation,
has been considered in
 Section 14.2, 14-16 to 14-21, for the explanation of the result in \cite{NPS} 
mentioned above.

First recall the probe method for Problem 1.

\subsection{The probe method}

The probe method formulated here is taken from \cite{INew} and should be called the Side A of the probe method
\cite{IHokkai}.
Given $x\in\Omega$ let $N_x$ denote the set of all non-self intersecting  piecewise linear curves $\sigma$ connecting
a point on $\partial\Omega$ and $x$ such that other points on $\sigma$ are in $\Omega$.  Each member in $N_x$ is called a {\it needle}
with a tip at $x$.   Let $G(y)$ denote the standard fundamental solution\footnote{There is a general version of the probe method, in which $G(y)$ is replaced with any solution
of the Laplace equation in $\Bbb R^3\setminus\{(0,0,0)\}$ with the property
$\int_V\,\vert\nabla G(y)\vert^2\,dy=\infty$ for any finite cone $V$ with vertex at $(0,0,0)$.  
See \cite{IProbeCarleman} and for other choices,  Proposition 4 in \cite{IProbe} for the stationary Schr\"odinger equation
and also \cite{IProbe2} for Navier's equation in the linear theory of elasticity.}
of the Laplace equation, that is
$$\begin{array}{ll}
\displaystyle
G(y)=\frac{1}{4\pi\vert y\vert}, & y\not=(0,0,0).
\end{array}
$$
A sequence of solutions  $\{v_n(\,\cdot\,;x)\}$ of the Laplace equation in $\Omega$ is called a {\it needle sequence}
for needle $\sigma\in N_x$ if  each member $v_n(\,\cdot\,;x)$ belongs to $H^1(\Omega)$ and 
the sequence $\{v_n(\,\cdot\,;x)\}$
converges to $G(\,\cdot\,-x)$ in the topology of $H^1_{\text{loc}}(\Omega\setminus\sigma)$
\footnote{The singular sources method never make use of this topology.  For each $x$ they fix a ``probing domain''
 $B_x$ in such a way that $x\notin\overline{B_x}$, $D\subset B_x$ and $\Omega\setminus\overline{B_x}$ is connected.
Instead of  $H^1_{\text{loc}}(\Omega\setminus\sigma)$-topology they require
$v_n(\,\cdot\,;x)\approx G(\,\cdot\,-x)$ in $H^1(B_x)$.  Theoretically our way is simple.  
However, I understand
that, it is impossible to have such a {\it sequence} numerically in advance, with respect to
this topology, instead, they considered a practical computation procedure of a {\it projection} onto $H^1(B_x)$ of the needle sequence.}.
It is well known that the existence of a needle sequence  for an arbitrary given $\sigma\in N_x$ has been ensured by applying the Runge approximation property for the Laplace equation.
We denote by the symbol $<\,\,,\,>$ the dual pairing between $H^{-\frac{1}{2}}(\partial\Omega)$ and $H^{\frac{1}{2}}(\partial\Omega)$.

The core part of  the Side A of  the probe method applied to Problem is stated as follows.

\proclaim{\noindent Theorem 1.}   

\noindent
(a)   Given $x\in\Omega\setminus\overline{D}$ and $\sigma\in N_x$, let $\{v_n(\,\cdot\,;x)\}$ be an arbitrary needle sequence for $\sigma$.  If $\sigma\cap\overline{D}=\emptyset$, then
we have
$$\displaystyle
\lim_{n\rightarrow\infty}<(\Lambda_0-\Lambda_D)(v_n(\,\cdot\,;x)\vert_{\partial\Omega}),v_n(\,\cdot\,;x)\vert_{\partial\Omega}>
=I(x),
\tag {1.2}
$$
where
$$\displaystyle
I(x)=\int_{\Omega\setminus\overline{D}}\,\vert\nabla w_x(y)\vert^2\,dy+\int_D\vert\nabla G(y-x)\vert^2\,dy
\tag {1.3}
$$
and the function $w_x=w\in H^1(\Omega\setminus\overline{D})$ is the weak solution of
$$\left\{\begin{array}{ll}
\displaystyle
\Delta w=0, & y\in\Omega\setminus\overline{D},
\\
\\
\displaystyle
\frac{\partial w}{\partial\nu}=-\frac{\partial}{\partial\nu}G(y-x), & y\in\partial D,
\\
\\
\displaystyle
w=0, & y\in\partial\Omega.
\end{array}
\right.
\tag {1.4}
$$

\noindent
(b)  We have, for each $\epsilon>0$
$$\displaystyle
\sup_{x\in\Omega\setminus\overline{D},\,\text{dist}\,(x,\partial D)>\epsilon}\,I(x)<\infty.
\tag {1.5}
$$

\noindent
(c) We have, for each $a\in\partial D$
$$\displaystyle
\lim_{x\rightarrow a}\,I(x)=\infty.
\tag {1.6}
$$

\endproclaim

The sequence $\{<(\Lambda_0-\Lambda_D)(v_n(\,\cdot\,;x)\vert_{\partial\Omega}), v_n(\,\cdot\,;x)\vert_{\partial\Omega}>\}$ in Theorem 1
is called {\it the indicator sequence} (\cite{INew}).
The formula (1.2) gives a way of calculating the value of the indicator function at an arbitrary point in $\Omega\setminus\overline{D}$.

Formula (1.2) is the direct consequence of two facts.

\noindent
(i)  The first is the well-posedness of the boundary value problem:
$$
\left\{
\begin{array}{ll}
\displaystyle
\Delta p=F, & x\in\Omega\setminus\overline{D},
\\
\\
\displaystyle
\frac{\partial p}{\partial\nu}=0, & x\in\partial D,
\\
\\
\displaystyle
p=0, & x\in\partial\Omega.
\end{array}
\right.
$$
That is, given an arbitrary $F\in L^2(\Omega\setminus\overline{D})$ the weak solution $p\in H^1(\Omega\setminus\overline{D})$ exists, unique and satisfies
$$\displaystyle
\Vert p\Vert_{L^2(\Omega\setminus\overline{D})}\le C\Vert F\Vert_{L^2(\Omega\setminus\overline{D})},
$$
where $C$ is a positive constant and independent of $F$.
Note that this estimate implies automatically an $L^2$ to $H^1$- estimate:
$$\displaystyle
\Vert p\Vert_{H^1(\Omega\setminus\overline{D})}\le C'\Vert F\Vert_{L^2(\Omega\setminus\overline{D})},
$$
where $C'$ is also a positive constant and independent of $F$.

\noindent
(ii)  The second is the well known decomposition formula \cite{E01}
$$\displaystyle
<(\Lambda_0-\Lambda_D)v\vert_{\partial\Omega},v\vert_{\partial\Omega}>
=\int_{\Omega\setminus\overline{D}}\vert\nabla(u-v)\vert^2\,dx+\int_D\vert\nabla v\vert^2\,dx,
$$
where $\Delta v=0$ in $\Omega$ and $u$ is the solution of (1.1) with $f=v$ on $\partial\Omega$.

From (i) and the convergence property of the needle sequence $v_n\rightarrow G(\,\cdot\,-x)$
in $H^1(D)$ under the assumption $\sigma\cap\overline{D}=\emptyset$, one has $u_n-v_n\rightarrow w_x$ in $H^1(\Omega\setminus\overline{D})$,
where $u_n$ solves (1.1) with $f=v_n$ on $\partial\Omega$ and thus (ii) yields (1.3).
The validity of (1.5) and (1.6) are easily checked since a combination of (i) and (ii) yields the estimates:
$$\displaystyle
\int_D\vert\nabla v\vert^2\,dx\le<(\Lambda_0-\Lambda_D)v\vert_{\partial\Omega},v\vert_{\partial\Omega}>
\le C''\int_D\vert\nabla v\vert^2\,dx,
\tag {1.7}
$$
where $C''$ is a positive constant independent of $v\in H^1(\Omega)$ satisfying $\Delta v=0$ in $\Omega$.
For these points, from the most recent view, see \cite{IR}.

The ``field'' $I(x)$ defined by (1.3) is called the {\it indicator function} for the probe method
 since by  (1.2) one can calculate the value of $I(x)$ from our observation data $\Lambda_0-\Lambda_D$;
by (1.5) $I(x)$ is  bounded  if $x$ is away from $\partial D$;
by (1.6) $I(x)$ blows up on $\partial D$.  

\noindent
The function $w_x$ is called the {\it reflected solution} by $D$.
  
$\quad$

{\bf\noindent Remark 1.1.}  Under the same assumption on $\sigma$ as (a),
we have also the formula:
$$\displaystyle
\lim_{m\rightarrow\infty}\lim_{n\rightarrow\infty}<(\Lambda_0-\Lambda_D)(v_n(\,\cdot\,;x)\vert_{\partial\Omega}),
v_m(\,\cdot\,;x)\vert_{\partial\Omega}>=I(x).
$$
This is proved as follows.
Let $u_n=u_n(\,\cdot\,;x)$ be the solution of (1.1) with $f=v_n(\,\cdot\,;x)$ on $\partial\Omega$.
We have
$$\begin{array}{l}
\displaystyle
\,\,\,\,\,\,
<(\Lambda_0-\Lambda_D)(v_n(\,\cdot\,;x)\vert_{\partial\Omega}),
v_m(\,\cdot\,;x)\vert_{\partial\Omega}>
\\
\\
\displaystyle
=\int_D\nabla v_n\cdot\nabla v_m\,dy+
\int_{\Omega\setminus\overline{D}}
\nabla (u_n-v_n)\cdot\nabla (u_m-v_m)\,dy.
\end{array}
$$
Since $v_n\rightarrow G(\,\cdot\,-x)$ in $H^1(D)$ and thus $u_n-v_n\rightarrow w_x$ in $H^1(\Omega\setminus\overline{D})$, we obtain the desired formula.

\subsection{The singular sources method and a Carleman function}

First let us consider what the singular sources method in a bounded domain should be.
For the purpose recall the essence of the original singular sources method for a sound-hard obstacle $D$ with $C^2$-boundary,
taken from \cite{P1,P3}.
The method employs the far field pattern of the scattered wave of an incident plane wave $e^{iky\cdot d}$, $d\in S^2$,  $k>0$ by the obstacle $D$, that is
$$\begin{array}{lll}
\displaystyle
w(r\vartheta;d,k)=\frac{e^{ikr}F(\vartheta, d;k)}{r}+O(\frac{1}{r^2}), & r=\vert x\vert\rightarrow\infty, & \vartheta\in S^2,
\end{array}
\tag {1.8}
$$
where $w=w(y;d,k)$ denotes the scattered wave and is characterized as the unique solution of
$$
\left\{
\begin{array}{ll}
(\Delta+k^2)w=0, & y\in\Bbb R^3\setminus\overline{D}\\
\\
\displaystyle
\frac{\partial w}{\partial\nu}=-\frac{\partial}{\partial\nu}e^{iky\cdot d}, & y\in\partial D,\\
\\
\displaystyle
r\left(\frac{\partial w}{\partial r}-ikw\,\right)\rightarrow 0, & r=\vert y\vert\rightarrow\infty.
\end{array}
\right.
$$
The coefficient $F(\vartheta,d;k)$ in the expansion (1.8) is called the far field pattern of $w$ at the direction $\vartheta$.
Hereafter $k$ is always fixed.
Next consider the scattered wave of the spherical wave generated by a point source located at a given point $x\in\Bbb R^3\setminus\overline{D}$.
It is the unique solution $\Phi_x$ of
$$
\left\{
\begin{array}{ll}
(\Delta+k^2)\Phi_x=0, & y\in\Bbb R^3\setminus\overline{D}\\
\\
\displaystyle
\frac{\partial \Phi_x}{\partial\nu}=-\frac{\partial}{\partial\nu}\Phi(y-x), & y\in\partial D,\\
\\
\displaystyle
r\left(\frac{\partial \Phi_x}{\partial r}-ik\Phi_x\,\right)\rightarrow 0, & r=\vert y\vert\rightarrow\infty,
\end{array}
\right.
\tag {1.9}
$$
where $\Phi$ denotes the outgoing fundamental solution of the Helmholtz equation given by
$$\begin{array}{ll}
\displaystyle
\Phi(y)=\frac{e^{ik\vert y\vert}}{4\pi\vert y\vert}, & y\not=(0,0,0).
\end{array}
$$
The theory of singular sources method consists of two facts.

$\quad$

{\bf\noindent Fact 1.}  One can calculate $\Phi_x(x)$ for $x$ close to $\partial D$ from the observation data formulated as $F(\vartheta,d;k)$ for all $d, \vartheta$.

$\quad$

{\bf\noindent  Fact 2.}  The $\Phi_x(x)$ blows up on $\partial D$, that is,
$$\displaystyle
\lim_{x\rightarrow a\in\partial D}\,\vert \Phi_x(x)\vert=\infty.
$$

$\quad$

The calculation processes in Fact 1 consist of three parts.   

\noindent
(a)  Calculation of  the near field $w(x;-d,k)$ for $x$ close to $\partial D$ and thus $e^{-ikx\cdot d}+w(x;-d,k)$
from the far field pattern $F(\vartheta,-d;k)$ for all $\vartheta$.

\noindent
(b)  Calculation of  the far field pattern of $\Phi_x(rd;k)$ as $r\rightarrow\infty$ by using  the {\it reciprocity formula} which connects
the far field pattern of  $\Phi_x$ with the near field $e^{-ikx\cdot d}+w(x;-d,k)$.

\noindent
(c)  Calculation of near field $\Phi_x(x)$ for $x$ close to $\partial D$ from the far field pattern of $\Phi_x(rd;k)$ as $r\rightarrow\infty$ calculated in (b).

The calculation processes in (a) and (c) are done, for example, by using the {\it point source method} \cite{P2}, or in principle, a classical far field pattern to near field continuation formula,
see \cite{CK}.   As a result one needs to take the limit twice.

In short, Fact 1 is a calculation process of the field $\Phi_x(x)$ from the observation data and Fact 2 means that the calculated field $\Phi_x(x)$ blows up on the surface of
the obstacle.  Comparing these facts with the statement in Theorem 1, we see that
the field $\Phi_x(x)$, $x\in\Bbb R^3\setminus\overline{D}$ plays the same role as the indicator function 
$I(x)$, $x\in\Omega\setminus\overline{D}$ in the probe method;
the far field pattern $F(\vartheta,d;k)$, $\vartheta,d\in S^2$ plays the role of
$\Lambda_0-\Lambda_D$.   We denote this by the following diagram.
$$\displaystyle
F(\vartheta,d;k),\vartheta,d\in S^2\rightarrow\text{The far-field pattern of $\Phi_x(y)$}\rightarrow\Phi_x(x)
$$

To transplant the idea of this singular sources method into the bounded domain case,
in this diagram, do the replacement $F(\vartheta,d;k),\vartheta,d\in S^2$ with $\Lambda_0-\Lambda_D$,
the far-field pattern of $\Phi_x(y)$ with the Cauch data of $X_x(y)$ on $\partial\Omega$, where
$X_x(y)$ is a field defined for $x\in\Omega\setminus\overline{D}$ to be determined later.  One gets the following diagram
$$\displaystyle
\Lambda_0-\Lambda_D\rightarrow\text{The Cauchy data of $X_x(y)$ on $\partial\Omega$}\rightarrow X_x(x)
$$
So, what is $X_x(y)$?  For this, suggested by Nakamura-Potthast-Sini \cite{NP} in the sound-soft obstacle case
and comparing the equations (1.4) with (1.9)
it is quite natural to choose
$X_x(y)=w_x(y)$, where $w_x(y)$ is the reflected solution of the probe method, which is the unique solution of (1.4).
Since $w_x=0$ on $\partial\Omega$ which should be the counterpart of  radiation condition
$\lim_{r\rightarrow\infty}r\left(\frac{\partial \Phi_x}{\partial r}-ik\Phi_x\,\right)=0$, knowing 
the Cauchy data of $w_x(y)$ on $\partial\Omega$ is equivalent
to knowing the Neumann data $\frac{\partial w_x}{\partial\nu}$ on $\partial\Omega$ and thus the Neumann data
shall play a role of the far field pattern of $\Phi_x(y)$.  For this recall the Rellich lemma \cite{CK}
and the uniqueness of the Cauchy problem for the Laplace equation.
Besides, by the interior regularity one can consider the value of $w_x(y)$ at $y=x\in\Omega\setminus\overline{D}$.
Thus one could obtain the following diagram.
$$\displaystyle
\Lambda_0-\Lambda_D\rightarrow\text{The Neumann data $\frac{\partial w_x}{\partial\nu}$ on $\partial\Omega$}\rightarrow w_x(x).
$$
The field $w_x(x)$ should be the indicator function in this expected singular sources method
if one could confirm Fact 2 for $\Phi_x(x)$ replaced with $w_x(x)$.

In fact, one gets the following result.

\proclaim{\noindent Theorem 2.}   Assume that  $\partial\Omega\in C^{1,1}$ and  $\partial D\in C^{1,1}$.

\noindent
(a)   Given $x\in\Omega\setminus\overline{D}$ and $\sigma\in N_x$ let $\{v_n(\,\cdot\,;x)\}$ be an arbitrary needle sequence for $\sigma$ with higher order regularity
$v_n(\,\cdot\,;x)\in H^2(\Omega)$.

If $\sigma\cap\overline{D}=\emptyset$, then
we have
$$\displaystyle
-\lim_{n\rightarrow\infty}
<(\Lambda_0-\Lambda_D)(v_n(\,\cdot\,;x)\vert_{\partial\Omega}),G_n(\,\cdot\,;x)\vert_{\partial\Omega}>
=w_x(x),
\tag {1.10}
$$
where $w_x(y)$ is the solution of (1.4) and $G_n(\,\cdot\,;x)$ is defined by
$$\begin{array}{ll}
\displaystyle
G_n(y;x)=G(y-x)-v_n(y;x), & y\in\Omega\setminus\{x\}.
\end{array}
\tag {1.11}
$$

\noindent
(b)  We have, for each $\epsilon>0$
$$\displaystyle
\sup_{x\in\Omega\setminus\overline{D},\,\text{dist}\,(x,\partial D)>\epsilon}\,\vert w_x(x)\vert<\infty.
\tag {1.12}
$$

\noindent
(c) We have, for each $a\in\partial D$
$$\displaystyle
\lim_{x\rightarrow a}\,w_x(x)=\infty.
\tag {1.13}
$$

\endproclaim

Comparing (1.10) with  (1.2), we clearly observe the difference of the calculation methods of  two indicator functions
$w_x(x)$ and $I(x)$.  However, the operation of taking the limit twice in the original singular sources method is avoided in 
formula (1.10).  This is an advantage of our approach.

By the well-posedness in $H^2(\Omega\setminus\overline{D})$ of (1.4) and the Sobolev imbedding \cite{Gr},
the validity of  (1.12) in (b) is clear.  For this, in Theorem 2, the higher regularity of $\partial\Omega$ and $\partial D$
is assumed.  In contrast to Theorem 2, Theorem 1 does not require such higher regularity.
Everything has been done in the framework of the weak solution.
Besides, in establishing (1.13) we never make use of an expression of the reflected solution $w_x(y)$ 
by solving an integral equation unlike the proof of Fact 2 in \cite{P1}, therein an expression 
of the scattered wave in terms of the potential theory \cite{CK} has been employed.

Some further remarks are in order.

$\quad$

{\bf\noindent Remark 1.2.}
Note that the calculation processes (a), (b) and (c) in Fact 1 (however, in (c) not closed to $\partial D$, instead, the surface of a domain 
that contains the obstacle) are also presented in \cite{I2}.  It is a reduction procedure for the application of the probe method.
Then the problem becomes the reconstruction problem of the obstacle from the data $\Phi_x(y)$, $x\not=y$ on the surface of a sphere
that contains the obstacle.  In \cite{Iwave} the reconstruction formula of the obstacle by applying the probe method to this reduced problem has been established.
For the {\it reciprocity formula} in (b)  see Lemma 1 on p.46 of \cite{R} for a sound-soft obstacle and also the {\it second method} in \cite{I2} for a sound-hard obstacle.

$\quad$

{\bf\noindent Remark 1.3.}  By the interior elliptic  regularity, we can always take the value of $w_x(y)$ at $y=x$.
By virtue of the regularity assumption on $\partial\Omega$ and $\partial D$, elliptic regularity up to boundary \cite{Gr} yields $w_x\in H^2(\Omega\setminus\overline{D})$. 
The validity of the formula (1.10) needs this higher order regularity.  This is contrast to formula (1.2).
And also one needs the higher order regularity for each member of needle sequence for the same purpose.
Such choice is always possible, by extending the original domain $\Omega$ larger and continuing needle $\sigma$ outside the extended domain.

$\quad$

{\bf\noindent Remark 1.4.}
Note that $G_n(\,\cdot\,;x)$ given by (1.11) satisfies
$$\begin{array}{ll}
\displaystyle
\Delta G_n(y;x)+\delta(y-x)=0, & y\in\Omega,
\end{array}
$$
and, for any compact set $K$ of $\Bbb R^3$ with $K\subset\Omega$ and $K\cap\sigma=\emptyset$ we have
$$\displaystyle
\Vert G_n(\,\cdot\,;x)\Vert_{H^1(K)}\rightarrow 0.
$$
Since $G_n(\,\cdot\,;x)$  is harmonic in $\Omega\setminus\sigma$, this convergence yields,
in particular\footnote{This is enough for our purpose} 
$$\displaystyle
\Vert G_n(\,\cdot\,;x)\Vert_{H^2(K)}\rightarrow 0.
$$
Thus if $K$ is chosen in such a way that $\overline D\subset K$ and $K\cap\sigma=\emptyset$, the function $G_n(\,\cdot\,;x)$ 
becomes a Carleman function (\cite{LRS,L}) for the domain $\Omega\setminus\overline{D}$ and the surface $\partial\Omega$
that gives us a direct calculation formula
of a solution
of the Laplace equation in $\Omega\setminus\overline{D}$ using only the Cauchy data on $\partial\Omega$.

$\quad$

{\bf\noindent Remark 1.5.}
Theorem 2 says that the field $w_x(x)$ should be called the indicator function in the singular sources method, which is described by using the notion of the probe method.

$\quad$

{\bf\noindent Remark 1.6.}
In \cite{NPS}, lines 5-8, on page 1551, there are some descriptions about how to calculate their indicator function from the Cauchy data.  They suggest to use the point source method (e.g., \cite{P2}), however,
no explicit calculation formula like (1.10)  is  given.

$\quad$

The regular part of a special Carleman function yields a concrete needle sequence.
This fact has been pointed out in \cite{IProbeCarleman}.
In this paper we employ the needle sequence as the {\it regular part} of a Carleman function.  
Formula  (1.10) can be considered as an application
of  this idea to the singular sources method.  This is also a different point from \cite{NPS}.

From Theorems 1 and 2 one can calculate {\it independently} both indicator functions in the probe and singular sources methods from the data $\Lambda_0-\Lambda_D$.

\subsection{The third indicator function}

Here we describe our main result.
Let $x\in\Omega\setminus\overline{D}$.
We start with a solution $W=W_x(y)$ of 
$$
\left\{
\begin{array}{ll}
\displaystyle
\Delta W=0, & y\in\Omega\setminus\overline{D},
\\
\\
\displaystyle
\frac{\partial W}{\partial\nu}=-\frac{\partial}{\partial\nu}G(y-x), & y\in\partial D,
\\
\\
\displaystyle
W=G(y-x), & y\in\partial\Omega.
\end{array}
\right.
\tag {1.14}
$$
Under the regularity assumption on $\partial\Omega$ and $\partial D$, we have $W_x\in H^2(\Omega\setminus\overline{D})$ and, by the interior regularity one can consider the value of $W_x(y)$
at $y=x$.

It is clear that $W_x(y)$ has the {\it trivial} and natural decomposition
$$\begin{array}{ll}
\displaystyle
W_x(y) 
& = \displaystyle w_x(y)+w_x^1(y),
\end{array}
\tag {1.15}
$$
where $w_x(y)$ is the same as Theorem 1 and $w_x^1(y)=w$ satisfies
$$
\left\{
\begin{array}{ll}
\displaystyle
\Delta w=0, & y\in\Omega\setminus\overline{D},
\\
\\
\displaystyle
\frac{\partial w}{\partial\nu}=0, & y\in\partial D,
\\
\\
\displaystyle
w=G(y-x), & y\in\partial\Omega.
\end{array}
\right.
\tag {1.16}
$$
In particular, (1.15) gives
$$\begin{array}{ll}
\displaystyle
W_x(x) 
& = \displaystyle w_x(x)+w_x^1(x).
\end{array}
\tag {1.17}
$$
In the following theorem we state this $W_x(x)$ has another decomposition.

\proclaim{\noindent Theorem 3.}
One can decompose $W_x(x)$ into another way:
$$\begin{array}{ll}
\displaystyle
W_x(x)
&
\displaystyle
=I(x)+I^1(x),
\end{array}
\tag {1.18}
$$
where  $I^1(x)$ is given by
$$\displaystyle
I^1(x)=\int_{\Omega\setminus\overline{D}}\,\vert\nabla w_x^1(z)\vert^2\,dz
+\int_{\Bbb R^3\setminus\overline\Omega}\,\vert\nabla G(z-x)\vert^2\,dz.
\tag {1.19}
$$

\endproclaim

$\quad$

{\bf\noindent Remark 1.7.}
From (1.3) and (1.19), one can rewrite (1.18):
$$\begin{array}{ll}
\displaystyle
W_x(x)
&
\displaystyle
=
\Vert\nabla w_x\Vert_{L^2(\Omega\setminus\overline{D})}^2
+\Vert\nabla w_x^1\Vert_{L^2(\Omega\setminus\overline{D})}^2
\\
\\
\displaystyle
&
\displaystyle
\,\,\,
+\Vert\nabla G(\,\cdot\,-x)\Vert_{L^2(D)}^2+\Vert\nabla G(\,\cdot\,-x)\Vert_{L^2(\Bbb R^3\setminus\overline{\Omega})}^2.
\end{array}
$$
Besides, it follows from the equations  (1.4) and (1.16) that
$$\displaystyle
\int_{\Omega\setminus\overline{D}}\,\nabla w_x^1(y)\cdot\nabla w_x(y)\,dy=0.
$$
Thus from (1.15) one obtains the decomposition
$$\displaystyle
W_x(x)=\Vert\nabla W_x\Vert_{L^2(\Omega\setminus\overline{D})}^2
+\Vert\nabla G(\,\cdot\,-x)\Vert_{L^2(\Bbb R^3\setminus
(\Omega\setminus\overline{D}))}^2,
$$
which corresponds to (1.3). This is the energy integral expression of the pointwise value of $W_x(y)$ on $y=x$.

The validity of the decomposition (1.18) is not a trivial fact at all.  The choice of $w_x^1$ corresponding
to $w_x$ is essential.
The proof employs a suitable needle sequence for each $x\in\Omega\setminus\overline{D}$ 
since the expressions (1.29) and (2.11) of $w_x$ and $w_x^1$ by using the needle sequence are used.
However, after having fixed $w_x^1$,  the proof of (1.18) itself can be done without making use of the needle sequence.
In Appendix we present an alternative direct proof of (1.18) which is independent of
the existence of the needle sequence.

\subsection{Consequences of Theorem 3}

Here we show how the statement (c)  of Theorem 2 follows from Theorem 3.

\noindent
The following facts can be easily checked by the well-posedness in $H^2(\Omega\setminus\overline{D})$  of  (1.16):
for each $\epsilon>0$
$$
\displaystyle
\sup_{x\in\Omega\setminus\overline{D},\, \text{dist}(x,\partial\Omega)>\epsilon} 
\vert w_x^1(x)\vert<\infty.
\tag {1.20}
$$
From (1.19), one gets
$$
\displaystyle
\sup_{x\in\Omega\setminus\overline{D},\, \text{dist}(x,\partial\Omega)>\epsilon} 
 I^1(x)<\infty.
 \tag {1.21}
$$

Since we have two-ways decompositions (1.17) and (1.18)  of $W_x(x)$, one can write
$$\displaystyle
w_x(x)-I(x)=I^1(x)-w_x^1(x).
\tag {1.22}
$$
Then, applying (1.20) and (1.21) to the right-hand side on (1.22), 
one gets,  each $\epsilon>0$
$$
\displaystyle
\sup_{x\in\Omega\setminus\overline{D},\, \text{dist}(x,\partial\Omega)>\epsilon} \vert w_x(x)-I(x)\vert<\infty.
\tag {1.23}
$$
Now (1.23) together with (1.6) implies (1.13).

Besides, let us consider the behaviour of $w_x^1(x)$ as $x\in\Omega\setminus\overline{D}$ approaches
a point on $\partial\Omega$.
Rewrite (1.22) as
$$\displaystyle
w_x^1(x)-I^1(x)=I(x)-w_x(x).
\tag {1.24}
$$
Applying (1.5) and (1.12) to the right-hand side on (1.24), one gets
$$
\displaystyle
\sup_{x\in\Omega\setminus\overline{D},\, \text{dist}(x,\partial D)>\epsilon} \vert w_x^1(x)-I^1(x)\vert<\infty.
\tag {1.25}
$$
Applying Fatou's lemma to the second term of the right-hand side on (1.19), we  have immediately
$$\displaystyle
\lim_{x\rightarrow b\in\partial\Omega}I^1(x)=\infty.
\tag {1.26}
$$
Thus this together (1.25) yields
$$
\displaystyle
\lim_{x\rightarrow b\in\partial\Omega}w_x^1(x)=\infty.
\tag {1.27}
$$

From these one can derive the property of $W_x(x)$.

\proclaim{\noindent Corollary 1.}  

\noindent
(a)  We have, for each $\epsilon>0$
$$
\left\{
\begin{array}{l}
\displaystyle
\sup_{x\in\Omega\setminus\overline{D},\, \text{dist}(x,\partial D)>\epsilon} 
(\vert W_x(x)-I^1(x)\vert+\vert W_x(x)-w_x^1(x)\vert)<\infty,\\
\\
\displaystyle
\sup_{x\in\Omega\setminus\overline{D},\, \text{dist}(x,\partial\Omega)>\epsilon} 
(\vert W_x(x)-I(x)\vert+\vert W_x(x)-w_x(x)\vert)<\infty.\\
 \end{array}
 \right.
$$

\noindent
(b) We have
$$\left\{
\begin{array}{l}
\displaystyle
\lim_{x\rightarrow a\in\partial D}W_x(x)=\infty,
\\
\\
\displaystyle
\lim_{x\rightarrow b\in\partial\Omega}W_x(x)=\infty.
\end{array}
\right.
$$

\endproclaim

So now we have the clear meaning of two indicator functions $w_x(x)$ and $I(x)$
through two decompositions (1.17) and (1.18) of $W_x(x)$.

\subsection{Calculating the third indicator function}

Both $I^1(x)$ and $w_x^1(x)$ can be calculated from $\Lambda_D$ as below.

\proclaim{\noindent Theorem 4.}
Let $x\in\Omega\setminus\overline{D}$.

\noindent
(a)  We have
$$\begin{array}{ll}
\displaystyle
I^1(x)
&
\displaystyle
=<\Lambda_D(G(\,\cdot-x)\,\vert_{\partial\Omega})-
\frac{\partial}{\partial\nu} G(\,\cdot-x)\vert_{\partial\Omega}, G(\,\cdot-x)\,\vert_{\partial\Omega}>.
\end{array}
\tag {1.28}
$$

\noindent
(b) Assume that $\partial \Omega\in C^{1,1}$ and $\partial D\in C^{1,1}$.
Let $\sigma\in N_x$ and  $\{v_n(\,\cdot\,;x)\}$ be an arbitrary needle sequence for $\sigma$ with higher order regularity $v_n(\,\cdot\,;x)\in H^2(\Omega)$.
If $\sigma\cap\overline{D}=\emptyset$, then we have
$$\displaystyle
w_x^1(x)=I^1(x)+\lim_{n\rightarrow\infty}
<(\Lambda_0-\Lambda_D)( v_n(\,\cdot\,;x)\vert_{\partial\Omega}), G(\,\cdot\,-x)\vert_{\partial\Omega}>.
\tag {1.29}
$$

\endproclaim

Given two functions $A(x)$ and $B(x)$ of  $x\in\Omega\setminus\overline{D}$, 
we write $A(x)\sim_{\partial D}\, B(x)$ as $x\rightarrow\partial D$ if for each $\epsilon>0$, we have
$$\displaystyle
\sup_{x\in\Omega\setminus\overline{D},\,\text{dist}\,(x,\partial\Omega)>\epsilon}\,\vert A(x)-B(x)\vert<\infty.
$$
And we write $A(x)\sim_{\partial\Omega} B(x)$ as $x\rightarrow\partial\Omega$ if
for each $\epsilon>0$, we have
$$\displaystyle
\sup_{x\in\Omega\setminus\overline{D},\,\text{dist}\,(x,\partial D)>\epsilon}\,\vert A(x)-B(x)\vert<\infty.
$$

$\quad$

Thus as a result one can calculate $W_x(x)$ itself from $\Lambda_D$ (and $\Lambda_0$).  
And we have:

\noindent
(i)  as $x\rightarrow \partial D$
$$\displaystyle
W_x(x)\sim_{\partial D} I(x)\sim_{\partial D} w_x(x);
$$

\noindent
(ii)  as $x\rightarrow \partial\Omega$
$$\displaystyle
W_x(x)\sim_{\partial\Omega} I^1(x)\sim_{\partial\Omega} w_x^1(x).
$$

So, in this sense,  $W_x(x)$ is also an indicator function and
the {\it generator} of the indicator functions of the probe and singular sources methods.  
This is the new intrinsic relationship between them.

$\quad$

{\bf\noindent Remark 1.8.}
Formula (1.28) implies that $I^1(x)$ can be calculated from $\Lambda_D$ acting only
$G(\,\cdot\,-x)\vert_{\partial\Omega}$ for all $x\in\Omega$
and thus no need of using any needle sequence for the calculation.

$\quad$

{\bf\noindent Remark 1.9.}
Formula (1.28) and (1.29) tell us that both $I^1(x)$ and $w_x^1(x)$ can be calculated 
independently from $I(x)$ and $w_x(x)$.  And formulae (1.26) and (1.27) correspond to (1.6) and (1.13), respectively.

$\quad$

This paper is organized as follows.
In Section 2, the choice of the third indicator function $W_x(x)$ is explained by using the notion of the probe method.
Then the proofs of  (1.10), (1.18), (1.28) and (1.29) are given.  In Section 3 
the {\it lifting}  $I(x,y)$ and $I^1(x,y)$ of the indicator functions $I(x)$, $I^1(x)$, respectively,
in the sense $I(x,y)=I(x)$ and $I^1(x,y)=I^1(x)$ on $y=x$ are introduced
and their analytical continuation properties are clarified.  Besides, it is shown that the sum  $I(x,y)+I^1(x,y)$
gives the decomposition of the {\it symmetrization} of $W_x(y)$ with respect to the variables $x$ and $y$.
This is an extension of (1.18) to the lifting.  Section 4 is devoted to considering 
the Side B of the singular sources method.  First, some discussion on difficulty
of clarifying the asymptotic behaviour of the sequence of 
the left-hand side of formula (1.10) is given.
Second, to overcome the difficulty a reformulation of  the probe and singular sources methods 
are introduced. Then the indicator functions defined therein, denoted by $I^*(x)$ and $w_x^*(x)$ on which both reformulated methods based, completely coincide with each other.  As a byproduct, it turns out that the reformulated singular sources method has also the Side B of the probe method, which is a characterization of the unknown obstacle by means of the indicator sequences.
In section 5 it is shown that the blowing up property on $\partial D$ of the indicator function $w_x(x)$ is invariant
with respect to a monotone perturbation of the domain $\Omega$.  This yields automatically
all the indicator functions $I(x)$, $W_x(x)$ and $I^*(x)=w_x^*(x)$ have the same property.
It is based on a recent result in \cite{IR}.  In the final section some remarks are given.
Appendix is devoted to an alternative direct proof  of  (1.18) starting with the classical representation formulae
of $w_x(y)$ and $w_x^1(y)$ by using their Cauchy data on $\partial\Omega$ and $\partial D$.

\section{Discussion ``why $W_x(x)$?'' and the proofs of  (1.10), (1.18), (1.28) and (1.29)}

The introduction of $W_x$ in Subsection 1.4, in particular, the choice of the boundary conditions on (1.14)
is not trivial.  In this section, we will explain the reason for this and also give proofs of (1.10) and (1.18).
We make use of the notion of the probe method together with the well known property of the indicator 
function $I(x)$.

First of all,  forget all the equations (1.4), (1.14) and (1.16).  
Instead, consider a general setting.

Let $f_x(y)$, $g_x(y)$ be families of functions of $y\in\Omega\setminus\{x\}$ indexed with $x\in\Omega$
and for simplicity, they belong to $C^{\infty}(\Bbb R^3\setminus\{x\})$.

Let $x\in\Omega\setminus\overline{D}$.
We start with the solution $W=W_x(y)$ of the boundary value problem
$$
\left\{
\begin{array}{ll}
\displaystyle
\Delta W=0, & y\in\Omega\setminus\overline{D},
\\
\\
\displaystyle
\frac{\partial W}{\partial\nu}=-\frac{\partial}{\partial\nu}f_x(y), & y\in\partial D,
\\
\\
\displaystyle
W=g_x(y), & y\in\partial\Omega.
\end{array}
\right.
\tag {2.1}
$$
Solution $W_x(y)$ has the {\it natural} decomposition:
$$\displaystyle
W_x(y)=w_x(y)+w_x^1(y),
$$
where the functions $w_x(y)$ and $w_x^1(y)$ satisfy
$$
\left\{
\begin{array}{ll}
\displaystyle
\Delta w_x=0, & y\in\Omega\setminus\overline{D},
\\
\\
\displaystyle
\frac{\partial w_x}{\partial\nu}=-\frac{\partial}{\partial\nu}f_x(y), & y\in\partial D,
\\
\\
\displaystyle
w_x=0, & y\in\partial\Omega
\end{array}
\right.
\tag {2.2}
$$
and
$$
\left\{
\begin{array}{ll}
\displaystyle
\Delta w_x^1=0, & y\in\Omega\setminus\overline{D},
\\
\\
\displaystyle
\frac{\partial w_x^1}{\partial\nu}=0, & y\in\partial D,
\\
\\
\displaystyle
w_x^1=g_x(y), & y\in\partial\Omega.
\end{array}
\right.
\tag {2.3}
$$

\subsection{Calculating $w_x^1(x)$}

First let us consider $w_x^1(y)$.
By the governing equations (2.3) and the definition of the Dirichlet-to-Neumann map, we have the expression
$$\begin{array}{ll}
\displaystyle
\frac{\partial w_x^1}{\partial\nu}(z)=\Lambda_D(g_x\vert_{\partial\Omega})(z), & z\in\partial\Omega.
\end{array}
\tag {2.4}
$$
Since $w_x^1=g_x$ on $\partial\Omega$, this means that the  Cauchy data of $w_x^1$ on $\partial\Omega$ can be calculated from $\Lambda_D$.
Thus, in principle, the value of $w_x^1$ at $y\in\Omega\setminus\overline{D}$ can be calculated from the Cauchy data.
In particular, one knows the value of $w_x^1(y)$ at $y=x$.

For the calculation, here we use the notion of the probe method.  
For $x\in\Omega\setminus\overline{D}$ choose a needle $\sigma\in N_x$ with $\sigma\cap \overline{D}=\emptyset$
and the needle sequence $\{v_m(\,\cdot\,;x)\}$ for $(x,\sigma)$. 
We can assume that $v_m\in H^2(\Omega)$ by extending the needle in a larger domain
and choosing a needle sequence for the extended needle.
By elliptic regularity up to boundary \cite{Gr}, we have $w_x^1\in H^2(\Omega\setminus\overline{D})$ and
interior regularity of $w_x^1$ yields
$$\displaystyle
w_x^1(x)=
\lim_{m\rightarrow\infty}
\int_{\partial\Omega}
\left(\Lambda_D(g_x\vert_{\partial\Omega})\,G_m(z;x)-\frac{\partial}{\partial\nu}
G_m(z;x)g_x(z)\,\right)\,dS(z),
\tag {2.5}
$$
where $G_m$ is the same as (1.11) and expression (2.4) is used.  By Remark 1.4, we see that
(2.5) is a consequence of Green's theorem, see Lemma 1.5.3.7 in \cite{Gr}.
The formula is nothing but a formula of Carleman type \cite{LRS}.

Rewrite the integral of the right-hand side on (2.5) as
$$\begin{array}{l}
\displaystyle
\,\,\,\,\,\,
\int_{\partial\Omega}
\left(\Lambda_D(g_x\vert_{\partial\Omega})\,G_m(z;x)-\frac{\partial}{\partial\nu}
G_m(z;x)g_x(z)\,\right)\,dS(z)
\\
\\
\displaystyle
=<\Lambda_D(g_x\vert_{\partial\Omega}), G(\,\cdot-x)\vert_{\partial\Omega}>
-<\Lambda_D(g_x\vert_{\partial\Omega}), v_m(\,\cdot\,;x)\vert_{\partial\Omega}>
\\
\\
\displaystyle
\,\,\,
-<\frac{\partial}{\partial\nu}G(\,\cdot-x)\vert_{\partial\Omega},g_x\vert_{\partial\Omega}>
+<\Lambda_0(v_m(\,\cdot\,;x)\vert_{\partial\Omega}),g_x\vert_{\partial\Omega}>.
\end{array}
$$
Thus (2.5) becomes
$$\begin{array}{ll}
\displaystyle
w_x^1(x)
&
\displaystyle
=
\lim_{m\rightarrow\infty}
<(\Lambda_0-\Lambda_D)(v_m(\,\cdot\,;x)\vert_{\partial\Omega}),g_x\vert_{\partial\Omega}>
\\
\\
\displaystyle
&
\displaystyle
\,\,\,
+<\Lambda_D(G(\,\cdot\,-x)\vert_{\partial\Omega})-\frac{\partial}{\partial\nu}G(\,\cdot-x)\vert_{\partial\Omega},
g_x\vert_{\partial\Omega}>.
\end{array}
\tag {2.6}
$$
Note that here we made use of the {\it symmetry} of $\Lambda_D$.
Also note that $g_x$ in this expression is still undetermined.
It should be emphasized that the essential part of the validity of this formula is the {\it existence}
of the needle sequence for an arbitrary needle $\sigma\in N_x$ with $\sigma\cap\overline{D}=\emptyset$.

\subsection{Derivation of (1.10)}

Next consider $w_x$.  Using the Cauchy data of $w_x$ on $\partial\Omega$,
similarly to the derivation of (2.5) we obtain
$$\displaystyle
w_x(x)=\lim_{m\rightarrow\infty}
\int_{\partial\Omega}
\frac{\partial w_x}{\partial\nu}(z)\,G_m(z;x)\,dS(z).
\tag {2.7}
$$
At the present time  the Neumann data $\frac{\partial w_x}{\partial\nu}$ on $\partial \Omega$ 
has no relation to $\Lambda_D$.  To make a direct relationship between $\frac{\partial w_x}{\partial\nu}\vert_{\partial\Omega}$
and the indicator function in the probe method, in (2.2)
we choose
$$\displaystyle
f_x(y)=G(y-x).
\tag {2.8}
$$
Now this $w_x$ coincides with that of Theorem 2 and we have
$$\displaystyle
w_x=\lim_{n\rightarrow \infty}(u_n(\,\cdot\,;x)-v_n(\,\cdot\,;x))
$$
in the topology of  $H^1(\Omega\setminus\overline{D})$, where $u=u_n$ solves (1.1) with $f=v_n(\,\cdot\,;x)$
on $\partial\Omega$.
Then we have the well known expression in the probe method:
$$\displaystyle
\frac{\partial w_x}{\partial\nu}=-\lim_{n\rightarrow\infty}
(\Lambda_0-\Lambda_D)(v_n(\cdot\,;x)\vert_{\partial\Omega}),
\tag {2.9}
$$
where the convergence is with respect to the topology of $H^{-\frac{1}{2}}(\partial\Omega)$.
Thus substituting (2.9) into (2.7), one can rewrite
$$\begin{array}{l}
\displaystyle
\,\,\,\,\,\,
w_x(x)
\\
\\
\displaystyle
=-\lim_{m\rightarrow\infty}<\lim_{n\rightarrow\infty}(\Lambda_0-\Lambda_D)(v_n(\cdot;x)\vert_{\partial\Omega}),
G_m(\,\cdot\,;x)\vert_{\partial\Omega}>
\\
\\
\displaystyle
=-\lim_{m\rightarrow\infty}<\lim_{n\rightarrow\infty}(\Lambda_0-\Lambda_D)(v_n(\,\cdot\,;x)\vert_{\partial\Omega}),
G(\,\cdot\,-x)\vert_{\partial\Omega}-v_m(\,\cdot\,;x)\vert_{\partial\Omega}>
\\
\\
\displaystyle
=-\lim_{n\rightarrow\infty}<(\Lambda_0-\Lambda_D)(v_n(\,\cdot\,;x)\vert_{\partial\Omega}),
G(\,\cdot\,-x)\vert_{\partial\Omega}>
\\
\\
\displaystyle
\,\,\,
+\lim_{m\rightarrow\infty}\lim_{n\rightarrow\infty}<(\Lambda_0-\Lambda_D)(v_n(\,\cdot\,;x)\vert_{\partial\Omega}),
v_m(\,\cdot\,;x)\vert_{\partial\Omega}>.
\end{array}
\tag {2.10}
$$
By Remark 1.1, (2.10) becomes
$$\displaystyle
w_x(x)
=I(x)-\lim_{n\rightarrow\infty}<(\Lambda_0-\Lambda_D)(v_n(\,\cdot\,;x)\vert_{\partial\Omega}),
G(\,\cdot\,-x)\vert_{\partial\Omega}>.
\tag {2.11}
$$
This together with (1.2) yields (1.10).

\subsection{Fixing $w_x^1$ and thus $W_x$.  Proof of  (1.18), (1.28) and (1.29)}

Here we choose $g_x(y)$ in (2.3) as
$$\displaystyle
g_x(y)=G(y-x).
\tag {2.12}
$$
Then, (2.6) becomes
$$\begin{array}{ll}
\displaystyle
w_x^1(x)
&
\displaystyle
=
\lim_{m\rightarrow\infty}
<(\Lambda_0-\Lambda_D)(v_m(\,\cdot\,;x)\vert_{\partial\Omega}),G(\,\cdot-x)\vert_{\partial\Omega}>
\\
\\
\displaystyle
&
\displaystyle
\,\,\,
+<\Lambda_D(G(\,\cdot\,-x)\vert_{\partial\Omega})-\frac{\partial}{\partial\nu}G(\,\cdot-x)\vert_{\partial\Omega}
,G(\,\cdot\,-x)\vert_{\partial\Omega}>.
\end{array}
\tag {2.13}
$$
Here write
$$
\displaystyle
-<\frac{\partial}{\partial\nu}G(\,\cdot-x)\vert_{\partial\Omega},G(\,\cdot\,-x)\vert_{\partial\Omega}>
=
-\int_{\partial\Omega}\frac{\partial}{\partial\nu}G(z-x)G(z-x)\,dS(z).
$$
Integration by parts yields
$$\displaystyle
-\int_{\partial\Omega}\frac{\partial}{\partial\nu}G(z-x)G(z-x)\,dS(z)
=\int_{\Bbb R^3\setminus\overline\Omega}\,\vert\nabla G(z-x)\vert^2\,dz.
\tag {2.14}
$$
From (1.16) and using integration by parts one gets
$$\displaystyle
<\Lambda_D(G(\,\cdot-x)\,\vert_{\partial\Omega}), G(\,\cdot-x)\,\vert_{\partial\Omega}>
=\int_{\Omega\setminus\overline{D}}\,\vert\nabla w_x^1(z)\vert^2\,dz.
\tag {2.15}
$$
Thus from (2.14) and (2.15) one gets
$$\displaystyle
\begin{array}{l}
\displaystyle
\,\,\,\,\,\,
<\Lambda_D(G(\,\cdot\,-x)\vert_{\partial\Omega})-\frac{\partial}{\partial\nu}G(\,\cdot-x)\vert_{\partial\Omega}
,G(\,\cdot\,-x)\vert_{\partial\Omega}>\\
\\
\displaystyle
=\int_{\Omega\setminus\overline{D}}\,\vert\nabla w_x^1(z)\vert^2\,dz
+\int_{\Bbb R^3\setminus\overline\Omega}\,\vert\nabla G(z-x)\vert^2\,dz.
\end{array}
\tag {2.16}
$$
A combination of (1.19) and (2.16) yields (1.28).
(1.29) is a consequence of (1.28) and (2.13).
(1.18) is a direct consequence of (1.29) and  (2.11) after choosing $f_x$ and $g_x$ as (2.8) and (2.12), respectively
in the equations (2.1), (2.2) and (2.3).

$\quad$

{\bf\noindent Remark 2.1.}
From (2.9) we have
$$\displaystyle
-\lim_{n\rightarrow\infty}<(\Lambda_0-\Lambda_D)(v_n(\,\cdot\,;x)\vert_{\partial\Omega}),
G(\,\cdot\,-x)\vert_{\partial\Omega}>
=\int_{\partial\Omega}\,\frac{\partial w_x}{\partial\nu}(z)G(z-x)\,dS(z)
$$
and thus (2.11) becomes
$$\displaystyle
w_x(x)
=I(x)+\int_{\partial\Omega}\,\frac{\partial w_x}{\partial\nu}(z)G(z-x)\,dS(z).
\tag {2.17}
$$
This is also a relationship between the two indicator functions $w_x(x)$ and $I(x)$
of the singular sources method  and the probe method.
The fact corresponding to (2.17) 
has been stated in Nakamura-Potthast 
\cite{NP} for sound-soft case
on page 14-21.
What they stated there can be expressed in our problem as follows:
as $x\rightarrow a\in\partial D$ 
$$\displaystyle
\int_{\partial\Omega}\,\frac{\partial w_x}{\partial\nu}(z)G(z-x)\,dS(z)
=O(1).
\tag {2.18}
$$
Note that by Lemma 4.1 in \cite{IHokkai}, we have
$$
\displaystyle
\Vert\nabla w_x\Vert_{L^2(\Omega\setminus\overline{D})}
\displaystyle
\ge C\frac{\Vert \nabla G(\,\cdot-x)\Vert_{L^2(D)}^2}{\Vert G(\,\cdot\,-x)\Vert_{H^1(D)}},
$$
where $C$ is a positive constant independent of $x\in\Omega\setminus\overline{D}$.
Since $\lim_{x\rightarrow a\in\partial D}\Vert\nabla G(\,\cdot\,-x)\Vert_{L^2(D)}=\infty$
and $\lim_{x\rightarrow a\in\partial D}\Vert G(\,\cdot\,-x)\Vert_{L^2(D)}<\infty$,
we have
$$\lim_{x\rightarrow a\in\partial D}
\frac{\Vert \nabla G(\,\cdot-x)\Vert_{L^2(D)}^2}{\Vert G(\,\cdot\,-x)\Vert_{H^1(D)}}=\infty
$$
and thus
$$\displaystyle
\lim_{x\rightarrow a\in\partial D}\Vert\nabla w_x\Vert_{L^2(\Omega\setminus\overline{D})}
=\infty.
$$
Then Poincar\'e inequality yields $\Vert w_x\Vert_{H^1(\Omega\setminus\overline{D})}\rightarrow\infty$.
So showing (2.18) is not a simple matter.  

Unfortunately, the author cannot find a description of this essential part in their book \cite{NP}.
Therefore, for a comparison of our third indicator function approach, in Section 5, we present  
{\it our own proof} of  (2.18) (see Remark 5.1), which is based on a global $L^2$-estimate of $w_x$
derived from lemma 2.2  in the recent article \cite{IR}.

Now, from our new point of view, the validity of (2.18) is {\it clear}.
From (2.17) and (1.22), we have
$$\displaystyle
\int_{\partial\Omega}\,\frac{\partial w_x}{\partial\nu}(z)G(z-x)\,dS(z)
=I^1(x)-w_x^1(x).
$$
This together with estimates (1.20) and (1.21) yields the desired conclusion.
This is also an advantage of introducing the solution of  (1.16) 
which was not considered at all in \cite{NP}.
The function $w_x^1(x)$ is just like an {\it additional line} in Euclidean Geometry.

$\quad$

{\bf\noindent Remark 2.2.}
Formula (2.11) and the symmetry of $\Lambda_D$ yield
$$\displaystyle
w_x(x)
=I(x)-\lim_{n\rightarrow\infty}<(\Lambda_0-\Lambda_D)(G(\,\cdot\,-x)\vert_{\partial\Omega}),
v_n(\,\cdot\,;x)\vert_{\partial\Omega}>.
$$
Thus $w_x(x)-I(x)$ can be calculated by knowing the data $\Lambda_D(G(\,\cdot\,-x)\vert_{\partial\Omega})$
and numerical values of  $v_n(\,\cdot\,;x)$ over $\partial\Omega$.

\section{Natural and inner decompositions}

Here we further consider the meaning of  two ways decomposition of $W_x(x)$.
For the purpose we start with the natural decomposition (1.15):
$$\begin{array}{ll}
\displaystyle
W_x(y)=w_x(y)+w_x^1(y), & (x,y)\in (\Omega\setminus\overline{D})\times(\Omega\setminus\overline{D}).
\end{array}
\tag {3.1}
$$

\subsection{Lifting the indicator functions}

First we introduce a lifting of the indicator function $I(x)$.

$\quad$

{\bf\noindent Definition 3.1.}
Define the function $I(x,y)$ of $(x,y)\in (\Omega\setminus\overline{D})\times(\Omega\setminus\overline{D})$
by the formula
$$\displaystyle
I(x,y)=\int_{\Omega\setminus\overline{D}}
\nabla w_x(z)\cdot\nabla w_y(z)\,dz+\int_D\nabla G(z-x)\cdot\nabla G(z-y)\,dz.
$$

$\quad$

We have $I(x)=I(x,y)\vert_{y=x}$.  In this sense, $I(x,y)$ should be called the lifting of $I(x)$.
The $I(x,y)$ can be calculated from $\Lambda_0-\Lambda_D$ as follows.

Given $(x,y)\in (\Omega\setminus\overline{D})\times(\Omega\setminus\overline{D})$
choose needles $\sigma\in N_x$ and $\sigma'\in N_y$ with $\sigma\cap\overline{D}=\sigma'\cap\overline{D}=
\emptyset$ and let $\{v_n(\,\cdot\,;x)\}$ and $\{v_n(\,\cdot\,;y)\}$ be arbitrary needle sequences
for $(x,\sigma)$ and $(y,\sigma')$, respectively.
We have 
$$\begin{array}{l}
\displaystyle
\,\,\,\,\,\,
<(\Lambda_0-\Lambda_D)(v_n(\,\cdot\,;x)\vert_{\partial\Omega}),
v_n(\,\cdot\,;y)\vert_{\partial\Omega}>
\\
\\
\displaystyle
=\int_{\Omega\setminus\overline{D}}
\nabla w_n(z;x)\cdot\nabla w_n(z;y)\,dz+\int_D\nabla v_n(z;x)\cdot\nabla v_n(z;y)\,dz.
\end{array}
$$
where $w_n(z;x)=u_n(z;x)-v_n(z;x)$ and $w_n(z;y)=u_n(z;y)-v_n(z;y)$;
$u_n(\,\cdot\,;y)$ and $u_n(\,\cdot:x)$ are the solutions
of (1.1) with $f=v_n(\,\cdot\,;x)$ and $f=v_n(\,\cdot;y)$, respectively.
Taking limit $n\rightarrow\infty$, we obtain
$$\begin{array}{l}
\displaystyle
\,\,\,\,\,\,
\lim_{n\rightarrow\infty}
<(\Lambda_0-\Lambda_D)(v_n(\,\cdot\,;x)\vert_{\partial\Omega}),
v_n(\,\cdot\,;y)\vert_{\partial\Omega}>
\\
\\
\displaystyle
=\int_{\Omega\setminus\overline{D}}
\nabla w_x(z)\cdot\nabla w_y(z)\,dz+\int_D\nabla G(z-x)\cdot\nabla G(z-y)\,dz.
\end{array}
$$
Thus $I(x,y)$ has the computation formula from $\Lambda_0-\Lambda_D$:
$$\displaystyle
I(x,y)=\lim_{n\rightarrow\infty}
<(\Lambda_0-\Lambda_D)(v_n(\,\cdot\,;x)\vert_{\partial\Omega}),
v_n(\,\cdot\,;y)\vert_{\partial\Omega}>.
\tag {3.2}
$$
Clearly, we have the symmetry
$$\displaystyle
I(x,y)=I(y,x).
$$
It follows from the governing equations for $w_x$ and $w_y$ that 
$$\begin{array}{ll}
\displaystyle
\int_{\Omega\setminus\overline{D}}
\nabla w_x(z)\cdot\nabla w_y(z)\,dz
&
\displaystyle
=\int_{\partial D}\frac{\partial}{\partial\nu}G(z-y) w_x(z)\,dS(z)
\\
\\
\displaystyle
&
\displaystyle
=\int_{\partial D}\frac{\partial}{\partial\nu}G(z-x) w_y(z)\,dS(z).
\end{array}
$$
Thus one gets another representation formula of $I(x,y)$, that is
$$\begin{array}{ll}\displaystyle
I(x,y)
&
\displaystyle
=\int_{\partial D}\frac{\partial}{\partial\nu}G(z-y) w_x(z)\,dS(z)
+\int_D\nabla G(z-x)\cdot\nabla G(z-y)\,dz
\\
\\
\displaystyle
&
\displaystyle
=\int_{\partial D}\frac{\partial}{\partial\nu}G(z-x) w_y(z)\,dS(z)
+\int_D\nabla G(z-x)\cdot\nabla G(z-y)\,dz.
\end{array}
\tag {3.3}
$$
Let $\Delta_x$ and $\Delta_y$ denote the Laplacian with respect to the variables 
$x=(x_1,x_2,x_3)$ and $y=(y_1,y_2,y_3)$, respectively.
From (3.3)  we  conclude that, for each fixed $x\in\Omega\setminus\overline{D}$
the function $I(x,\,\cdot\,)\in C^{\infty}(\Omega\setminus\overline{D})$ and satisfies
$$\begin{array}{ll}
\displaystyle
\Delta_y I(x,y)=0, & y\in\Omega\setminus\overline{D}
\end{array}
\tag {3.4}
$$
and for each fixed $y\in\Omega\setminus\overline{D}$
the function $I(\,\cdot\,,y)\in C^{\infty}(\Omega\setminus\overline{D})$ and satisfies
$$\begin{array}{ll}
\displaystyle
\Delta_x I(x,y)=0, & x\in\Omega\setminus\overline{D}.
\end{array}
\tag {3.5}
$$
Thus from (3.4) and (3.5) together with the weak unique continuation theorem for the Laplace equation 
we obtain

\proclaim{\noindent Proposition 3.1.}
Let $U$ and $V$ be non empty open subsets of  $\Omega\setminus\overline{D}$.
Then the values of $I(x,y)$, $(x,y)\in U\times V$ uniquely determine all the values
of $I(x,y)$, $(x,y)\in(\Omega\setminus\overline{D})^2$
and thus the indicator function $I(x)=I(x,x)$, $x\in\Omega\setminus\overline{D}$.
\endproclaim

Next we introduce a lifting of $I^1(x)$.

$\quad$

{\bf\noindent Definition 3.2.}
Define the function $I^1(x,y)$ of  $(x,y)\in (\Omega\setminus\overline{D})\times(\Omega\setminus\overline{D})$
by the formula
$$\displaystyle
I^1(x,y)=\int_{\Omega\setminus\overline{D}}\,\nabla w_x^1(z)\cdot\nabla w_y^1(z)\,dz
+\int_{\Bbb R^3\setminus\overline\Omega}\,\nabla G(z-x)\cdot\nabla G(z-y)\,dz.
$$

$\quad$

Clearly we have the symmetry
$$\displaystyle
I^1(x,y)=I^1(y,x).
$$

Similarly to (2.15), using the governing equations of $w_x$ and $w_y$, one has
the expression
$$\displaystyle
\int_{\Omega\setminus\overline{D}}\,\nabla w_x^1(z)\cdot\nabla w_y^1(z)\,dz
=<\Lambda_D(G(\,\cdot-x)\,\vert_{\partial\Omega}), G(\,\cdot-y)\,\vert_{\partial\Omega}>
$$
and also similarly to (2.14) we have
$$\displaystyle
-<\frac{\partial}{\partial\nu} G(\,\cdot-x)\vert_{\partial\Omega},G(\,\cdot\,-y)\vert_{\partial\Omega}>
=\int_{\Bbb R^3\setminus\overline{\Omega}}\,\nabla G(z-x)\cdot\nabla G(z-y)\,dz.
$$
Thus $I^1(x,y)$ becomes
$$\begin{array}{ll}
\displaystyle
I^1(x,y)
&
\displaystyle
=<\Lambda_D(G(\,\cdot-x)\,\vert_{\partial\Omega})-\frac{\partial}{\partial\nu} G(\,\cdot-x)\vert_{\partial\Omega}
, G(\,\cdot-y)\,\vert_{\partial\Omega}>
\\
\\
\displaystyle
&
\displaystyle
=<\Lambda_D(G(\,\cdot-y)\,\vert_{\partial\Omega})
-\frac{\partial}{\partial\nu} G(\,\cdot-y)\vert_{\partial\Omega}, G(\,\cdot-x)\,\vert_{\partial\Omega}>.
\end{array}
\tag {3.6}
$$
Therefore, it is easy to see that, for each fixed $x\in\Omega\setminus\overline{D}$
the function $I^1(x,\,\cdot\,)\in C^{\infty}(\Omega\setminus\overline{D})$ and satisfies
$$\begin{array}{ll}
\displaystyle
\Delta_y I^1(x,y)=0, & y\in\Omega\setminus\overline{D},
\end{array}
\tag {3.7}
$$
and for each fixed  $y\in\Omega\setminus\overline{D}$
the function $I^1(\,\cdot\,,x)\in C^{\infty}(\Omega\setminus\overline{D})$ and satisfies
$$\begin{array}{ll}
\displaystyle
\Delta_x I^1(x,y)=0, & x\in\Omega\setminus\overline{D}.
\end{array}
\tag {3.8}
$$
Note that besides, the expression (3.6) gives the harmonic continuation of  $I^1(x,\,\cdot\,)$
and  $I^1(\,\cdot\,,y)$ onto $\Omega$ for each fixed $x\in\Omega\setminus\overline{D}$
and $y\in\Omega\setminus\overline{D}$, respectively\footnote{However, in this paper we do not make use of  this fact
anywhere.}.

From (3.7), (3.8) and the weak unique continuation theorem for the Laplace equation we obtain

\proclaim{\noindent Proposition 3.2.}
Let $U$ and $V$ be non empty open subsets of  $\Omega\setminus\overline{D}$.
Then the values of $I^1(x,y)$, $(x,y)\in U\times V$ uniquely determine all the values
of $I^1(x,y)$, $(x,y)\in(\Omega\setminus\overline{D})^2$
and thus the indicator function $I^1(x)=I^1(x,x)$, $x\in\Omega\setminus\overline{D}$.
\endproclaim

\subsection{Inner decomposition}

In this subsection we claify the relationship between $I(x,y)$, $I^1(x,y)$, $w_x(y)$, $w_x^1(y)$ and $W_x(y)$.

\proclaim{\noindent Theorem 5.}
Given $(x,y)\in (\Omega\setminus\overline{D})\times(\Omega\setminus\overline{D})$
choose needles $\sigma\in N_x$ and $\sigma'\in N_y$ with $\sigma\cap\overline{D}=\sigma'\cap\overline{D}=
\emptyset$ and let $\{v_n(\,\cdot\,;x)\}$ and $\{v_n(\,\cdot\,;y)\}$ be arbitrary needle sequences
for $(x,\sigma)$ and $(y,\sigma')$, respectively.
We have
$$\displaystyle
w_x(y)=I(x,y)-\lim_{n\rightarrow\infty}<(\Lambda_0-\Lambda_D)(v_n(\,\cdot\,;x)\vert_{\partial\Omega}),
G(\,\cdot\,-y)\vert_{\partial\Omega}>,
\tag {3.9}
$$
$$\displaystyle
w_x^1(y)
=I^1(x,y)+\lim_{n\rightarrow\infty}
<(\Lambda_0-\Lambda_D)(v_n(\,\cdot\,;y)\vert_{\partial\Omega}),G(\,\cdot\,-x)\vert_{\partial\Omega}>
\tag {3.10}
$$
and 
$$\begin{array}{ll}
\displaystyle
W_x(y)
&
\displaystyle
=I(x,y)+I^1(x,y)
\\
\\
\displaystyle
&
\displaystyle
\,\,\,
+\lim_{n\rightarrow\infty}
<(\Lambda_0-\Lambda_D)(v_n(\,\cdot\,;y)\vert_{\partial\Omega}),G(\,\cdot\,-x)\vert_{\partial\Omega}>
\\
\\
\displaystyle
&
\displaystyle
\,\,\,
-\lim_{n\rightarrow\infty}<(\Lambda_0-\Lambda_D)(v_n(\,\cdot\,;x)\vert_{\partial\Omega}),
G(\,\cdot\,-y)\vert_{\partial\Omega}>.
\end{array}
\tag {3.11}
$$
Therefore the term $I(x,y)+I^1(x,y)$ coincides with the symmetric part of $W_x(y)$ with respect to $x$ and $y$,
that is, 
$$\displaystyle
\frac{W_x(y)+W_y(x)}{2}=I(x,y)+I^1(x,y).
\tag {3.12}
$$

\endproclaim

{\it\noindent Proof.}
We have
$$\displaystyle
w_x(y)=\lim_{m\rightarrow\infty}
\int_{\partial\Omega}
\frac{\partial w_x}{\partial\nu}(z)\,G_m(z;y)\,dS(z)
\tag {3.13}
$$
and also the well known expression in the probe method
$$\displaystyle
\frac{\partial w_x}{\partial\nu}(z)=-\lim_{n\rightarrow\infty}
(\Lambda_0-\Lambda_D)(v_n(\cdot\,;x)\vert_{\partial\Omega})(z),
\tag {3.14}
$$
where the convergence is with respect to the topology of $H^{-\frac{1}{2}}(\partial\Omega)$.
Thus substituting (3.14) into (3.13), one can rewrite
$$\begin{array}{l}
\displaystyle
\,\,\,\,\,\,
w_x(y)
\\
\\
\displaystyle
=-\lim_{m\rightarrow\infty}<\lim_{n\rightarrow\infty}(\Lambda_0-\Lambda_D)(v_n(\cdot;x)\vert_{\partial\Omega}),
G_m(\,\cdot\,;y)\vert_{\partial\Omega}>
\\
\\
\displaystyle
=-\lim_{m\rightarrow\infty}<\lim_{n\rightarrow\infty}(\Lambda_0-\Lambda_D)(v_n(\,\cdot\,;x)\vert_{\partial\Omega}),
G(\,\cdot\,-y)\vert_{\partial\Omega}-v_m(\,\cdot\,;y)\vert_{\partial\Omega}>
\\
\\
\displaystyle
=-\lim_{n\rightarrow\infty}<(\Lambda_0-\Lambda_D)(v_n(\,\cdot\,;x)\vert_{\partial\Omega}),
G(\,\cdot\,-y)\vert_{\partial\Omega}>
\\
\\
\displaystyle
\,\,\,
+\lim_{m\rightarrow\infty}\lim_{n\rightarrow\infty}<(\Lambda_0-\Lambda_D)(v_n(\,\cdot\,;x)\vert_{\partial\Omega}),
v_m(\,\cdot\,;y)\vert_{\partial\Omega}>.
\end{array}
\tag {3.15}
$$
Using  (3.2) together with a remark similar to Remark 1.1,  we see that (3.15) yields (3.9) which is a generalization of (2.11).

Next, recalling (2.5) with $g_x(z)$ given by (2.12), we have
$$\displaystyle
w_x^1(y)=
\lim_{m\rightarrow\infty}
\int_{\partial\Omega}
\left(\Lambda_D(G(\,\cdot\,-x)\vert_{\partial\Omega})\,G_m(z;y)-\frac{\partial}{\partial\nu}
G_m(z;y)G(z-x)\,\right)\,dS(z),
\tag {3.16}
$$
where $G_m(z;y)=G(z-y)-v_m(z;y)$.
Rewrite the integral of the right-hand side on (3.16), one gets
$$\begin{array}{l}
\displaystyle
\,\,\,\,\,\,
\int_{\partial\Omega}
\left(\Lambda_D(G(\,\cdot\,-x)\vert_{\partial\Omega})\,G_m(z;y)-\frac{\partial}{\partial\nu}
G_m(z;y)G(z-x)\,\right)\,dS(z)
\\
\\
\displaystyle
=<\Lambda_D(G(\,\cdot\,-x)\vert_{\partial\Omega}), G(\,\cdot-y)\vert_{\partial\Omega}>
-<\Lambda_D(G(\,\cdot\,-x)\vert_{\partial\Omega}), v_m(\,\cdot\,;y)\vert_{\partial\Omega}>
\\
\\
\displaystyle
\,\,\,
-<\frac{\partial}{\partial\nu}G(\,\cdot-y)\vert_{\partial\Omega},G(\,\cdot\,-x)\vert_{\partial\Omega}>
+<\Lambda_0(v_m(\,\cdot\,;y)\vert_{\partial\Omega},G(\,\cdot\,-x)\vert_{\partial\Omega}>.
\end{array}
$$
Thus (3.16) becomes
$$\begin{array}{ll}
\displaystyle
w_x^1(y)
&
\displaystyle
=
\lim_{m\rightarrow\infty}
<(\Lambda_0-\Lambda_D)(v_m(\,\cdot\,;y)\vert_{\partial\Omega}),G(\,\cdot\,-x)\vert_{\partial\Omega}>
\\
\\
\displaystyle
&
\displaystyle
\,\,\,
+<\Lambda_D(G(\,\cdot\,-y)\vert_{\partial\Omega})-\frac{\partial}{\partial\nu}G(\,\cdot-y)\vert_{\partial\Omega},G(\,\cdot\,-x)\vert_{\partial\Omega}>.
\end{array}
\tag {3.17}
$$
Note that here we made use of the {\it symmetry} of $\Lambda_D$.
Then, substituting (3.6) into (3.17), we obtain (3.10).

Finally (3.1),  (3.9) and (3.10) give us the expression (3.11).

\noindent
$\Box$

$\quad$

{\bf\noindent Remark 3.1.}

\noindent
(a)  We call the decomposition (3.12) the {\it inner decomposition} of 
the symmetric part of $W_x(y)$
with respect to $x$ and $y$. 
Note that we have, for all $x,y\in\Omega\setminus\overline{D}$
$$\displaystyle
\int_{\Omega\setminus\overline{D}}\nabla w_x(z)\cdot\nabla w_y^1(z)\,dz=\int_{\Omega\setminus\overline{D}}\nabla w_y(z)\cdot\nabla w_x^1(z)\,dz=0.
$$
Thus one gets
$$\displaystyle
I(x,y)+I^1(x,y)=\int_{\Omega\setminus\overline{D}}\nabla W_x(z)\cdot\nabla W_y(z)\,dz
+\int_{\Bbb R^3\setminus(\Omega\setminus\overline{D})}\nabla G(z-x)\cdot\nabla G(z-y)\,dz.
$$
Then, (3.12) yields the nearly {\it self-contained} identity:
$$\displaystyle
\frac{W_x(y)+W_y(x)}{2}=\int_{\Omega\setminus\overline{D}}\nabla W_x(z)\cdot\nabla W_y(z)\,dz
+\int_{\Bbb R^3\setminus(\Omega\setminus\overline{D})}\nabla G(z-x)\cdot\nabla G(z-y)\,dz.
$$

\noindent
(b)  From (3.9) and (3.10)  and the symmetry of $I(x,y)$ and $I^1(x,y)$ we have
the {\it twisted decomposition}
$$\displaystyle
w_x^1(y)+w_y(x)=I(x,y)+I^1(x,y)
\tag {3.18}
$$
or equivalently
$$\displaystyle
w_y^1(x)+w_x(y)=I(x,y)+I^1(x,y).
\tag{3.19}
$$
Thus one gets
$$
\displaystyle
\Delta_x(w_x^1(y))=\Delta_x(w_x(y))=\Delta_x(W_x(y))=0.
\tag {3.20}
$$

The proof of  (3.20) in principle, may be done directly by using the governing equation for $W_x(y)$ or $w_x(y)$ and
$w_x^1(y)$, the well-posedness of the boundary value problem in some function spaces.
However, the proof mentioned above is more straightforward
since it is based on
the formulae (3.18) (and (3.19)), and  equations (3.4), (3.5), (3.7) and (3.8).
The equations (3.20) given here are based on the differentiation of the expressions (3.3) and (3.6) under integral symbol {\it point-wisely}.

It should be emphasized that the proof of (3.18) itself presented here is not trivial since it makes use of the Runge approximation property for the Laplace equation which ensures the existence of needle sequences.
However, as can be seen from an alternative proof of (1.18) in Appendix, see (A.3) and (A.7), 
one can make it avoided.
In this sense, one can say that the Runge approximation property is useful to find an integral identity,
however, after finding the identity, in some case, one can find a direct elementary proof and we should do it.

\noindent
(c)  Applying the unique continuation theorem for the Laplace equation, we obtain the following fact
which corresponds to Propositions 3.1 and 3.2.
\proclaim{\noindent Proposition 3.3.}
Let $U$ and $V$ be non empty open subsets of  $\Omega\setminus\overline{D}$.
Then the values of $w_x(y)$, $w_x^1(y)$ and $W_x(y)$, $(x,y)\in U\times V$ uniquely determine all the values
of $w_x(y)$, $w_x^1(y)$ and $W_x(y)$, $(x,y)\in(\Omega\setminus\overline{D})^2$
and thus those of  $w_x(x)$, $w_x^1(x)$ and the third indicator function $W_x(x)$, $x\in\Omega\setminus\overline{D}$, respectively.
\endproclaim

\noindent
Note that Propositions 3.1, 3.2 and 3.3 are just the statement on the uniqueness and any constructive
procedure is not given.
However, those are surely new facts never mentioned in the previous probe and singular sources methods.

\section{The Side B of the singular sources method}

\subsection{A discussion on the Side B}

The Side B of the probe method is concerned with the asymptotic 
bahaviour of the sequence 
$$\displaystyle
\{<(\Lambda_0-\Lambda_D)(v_n(\,\cdot\,;x)\vert_{\partial\Omega}),v_n(\,\cdot\,;x)\vert_{\partial\Omega}>\}
\tag {4.1}
$$ 
as $n\rightarrow\infty$ when $\sigma\cap\overline{D}\not=\emptyset$.
It is based on the system (1.7) and the blowing up property of the needle sequence on the needle which has been established in \cite{INew}.

In \cite{INew}, we have the following result.

\proclaim{\noindent Theorem 6.}  Let $x\in\Omega$ and $\sigma\in N_x$.
Assume that one of two cases (a) and (b) are satisfied:

\noindent
(a) $x\in\overline{D}$;

\noindent
(b) $x\in\Omega\setminus\overline{D}$ and $\sigma\cap D\not=\emptyset$.

\noindent
Then, for any needle sequence $\{v_n(\,\cdot\,;x)\}$ for $(x,\sigma)$
we have 
$$\displaystyle
\lim_{n\rightarrow\infty}
<(\Lambda_0-\Lambda_D)(v_n(\,\cdot\,;x)\vert_{\partial\Omega}),v_n(\,\cdot\,;x)\vert_{\partial\Omega}>=\infty.
$$

\endproclaim

In particular, (a) of Theorem 6 yields, if $x\in\overline{D}$, then
for any needle $\sigma\in N_x$ and needle sequence $\{v_n(\,\cdot\,;x)\}$ for $(x,\sigma)$
$$\displaystyle
\lim_{n\rightarrow\infty}
<(\Lambda_0-\Lambda_D)(v_n(\,\cdot\,;x)\vert_{\partial\Omega}),v_n(\,\cdot\,;x)\vert_{\partial\Omega}>=\infty.
$$
On the other hand, by (a) of Theorem 1, the connectedness of $\Omega\setminus\overline{D}$ and
the existence of the needle sequence, we know, if $x\in\Omega\setminus\overline{D}$, then
there exists a needle $\sigma\in N_x$ such that, for all needle sequences $\{v_n(\,\cdot\,;x)\}$ for $(x,\sigma)$
the limit
$$\displaystyle
\lim_{n\rightarrow\infty}
<(\Lambda_0-\Lambda_D)(v_n(\,\cdot\,;x)\vert_{\partial\Omega}),v_n(\,\cdot\,;x)\vert_{\partial\Omega}>
$$
exists.

As a direct corollary of  Theorems 1 and  6, in \cite{INew} the following characterization of the obstacle $D$ in terms of
the behaviour of the sequence (4.1) is given.

$\quad$

\proclaim{\noindent Corollary 2.}  A point $x\in\Omega$ belongs to the set $\Omega\setminus\overline{D}$ if and only if there exists a needle $\sigma\in N_x$ and a needle sequence $\{v_n(\,\cdot\,;x)\}$ for $(x,\sigma)$
such that
$$\displaystyle
\sup_n\,<(\Lambda_0-\Lambda_D)(v_n(\,\cdot\,;x)\vert_{\partial\Omega}),v_n(\,\cdot\,;x)\vert_{\partial\Omega}><\infty.
$$

\endproclaim

$\quad$

Or we have an explicit characterization of $\overline{D}$ itself. 

\proclaim{\noindent Corollary 3.}
A point $x\in\Omega$ belongs to the set $\overline{D}$ if and only if
for any needle $\sigma\in N_x$ and needle sequence $\{v_n(\,\cdot\,;x)\}$ for $(x,\sigma)$
we have
$$\displaystyle
\lim_{n\rightarrow\infty}<(\Lambda_0-\Lambda_D)(v_n(\,\cdot\,;x)\vert_{\partial\Omega}),v_n(\,\cdot\,;x)\vert_{\partial\Omega}>=\infty.
$$
\endproclaim
{\it\noindent Proof.}  The ``only if part''  is already explained.  Next assume that $x\not\in\overline{D}$.
Then as explained above, we have:
there exists a needle $\sigma\in N_x$ such that, for all needle sequences $\{v_n(\,\cdot\,;x)\}$ for $(x,\sigma)$
the limit
$$\displaystyle
\lim_{n\rightarrow\infty}
<(\Lambda_0-\Lambda_D)(v_n(\,\cdot\,;x)\vert_{\partial\Omega}),v_n(\,\cdot\,;x)\vert_{\partial\Omega}>
$$
exists.  Since, the needle sequence always {\it exists} for $(x,\sigma)$, one concludes:
there exists a needle $\sigma\in N_x$ and needle sequence $\{v_n(\,\cdot\,;x)\}$ for $(x,\sigma)$
such that the limit
$$\displaystyle
\lim_{n\rightarrow\infty}
<(\Lambda_0-\Lambda_D)(v_n(\,\cdot\,;x)\vert_{\partial\Omega}),v_n(\,\cdot\,;x)\vert_{\partial\Omega}>
$$
exists and thus the statement 
$$\displaystyle
\lim_{n\rightarrow\infty}<(\Lambda_0-\Lambda_D)(v_n(\,\cdot\,;x)\vert_{\partial\Omega}),v_n(\,\cdot\,;x)\vert_{\partial\Omega}>=\infty,
$$
is not valid.  This completes the proof of  the ``if  part''.

\noindent
$\Box$

\noindent
Note that the proof of ``if part'' presented above  is not a simple matter since it makes use of the existence
of the needle sequence.

Corollaries 2 and 3 clearly explain the reason why we call the sequence  (4.1)
the {\it indicator sequence} \cite{INew}.   Note that those corollaries yield automatically
the uniqueness of $D$.

In the singular sources method as formulated in Subsection 1.2, taking an account of formula (1.10),
the corresponding sequence should be the one given by
$$\displaystyle
\{-<(\Lambda_0-\Lambda_D)(v_n(\,\cdot\,;x)\vert_{\partial\Omega}),G_n(\,\cdot\,;x)\vert_{\partial\Omega}>\}.
\tag {4.2}
$$
Since we have
$$\begin{array}{l}
\displaystyle
\,\,\,\,\,\,
-<(\Lambda_0-\Lambda_D)(v_n(\,\cdot\,;x)\vert_{\partial\Omega}),G_n(\,\cdot\,;x)\vert_{\partial\Omega}>
\\
\\
\displaystyle
=<(\Lambda_0-\Lambda_D)(v_n(\,\cdot\,;x)\vert_{\partial\Omega}),v_n(\,\cdot\,;x)\vert_{\partial\Omega}>
\\
\\
\displaystyle
\,\,\,
-<(\Lambda_0-\Lambda_D)(v_n(\,\cdot\,;x)\vert_{\partial\Omega}),G(\,\cdot\,-x)\vert_{\partial\Omega}>,
\end{array}
$$
to translate the result in the Side B of the probe method into this, we have to show
the boundedness of the sequence (4.2) when $\sigma\cap\overline{D}\not=\emptyset$.

Integration by parts, the symmetry of $\Lambda_D$ and equations (1.16) give us 
$$\begin{array}{l}
\displaystyle
\,\,\,\,\,\,
<(\Lambda_0-\Lambda_D)(v_n(\,\cdot\,;x)\vert_{\partial\Omega}),G(\,\cdot\,-x)\vert_{\partial\Omega}>\\
\\
\displaystyle
=<\Lambda_0(v_n(\,\cdot\,;x)\vert_{\partial\Omega}, w_x^1\vert_{\partial\Omega}>
-<\Lambda_D(w_x^1\vert_{\partial\Omega}), v_n(\,\cdot\,;x)\vert_{\partial\Omega}>
\\
\\
\displaystyle
=\int_{\partial D}\,\frac{\partial}{\partial\nu}v_n(y;x)\,w_x^1(y)\,dS(y).
\end{array}
\tag {4.3}
$$
It seems that clarifying the asymptotic behaviour of this last integral as $n\rightarrow\infty$
is not clear when $\sigma\cap\overline{D}\not=\emptyset$.
The reason is:
the convergence   $v_n(\,\cdot\,;x)\rightarrow G(\,\cdot\,-x)$ in $H^1_{\text{loc}}(\Omega\setminus\sigma)$
does not control the behaviour of $v_n(\,\cdot\,;x)$ on $\partial D$ if $\sigma\cap\overline{D}\not=\emptyset$.
For example, consider the simplest situation where 
 $x\in\partial D$ and $(\sigma\setminus\{x\})\cap\overline{D}=\emptyset$.
 Then we have, for all $y\in\partial D\setminus\{x\}$
 $$\displaystyle
 \lim_{n\rightarrow\infty}\frac{\partial}{\partial\nu} v_n(y;x)=\frac{\partial }{\partial\nu}G(y-x).
 $$
However, the original convergence of needle sequence in $H^1_{\text{loc}}(\Omega\setminus\sigma)$
does not ensure the uniform convergence in a neighbourhood of $x$.
Besides we cannot expect the existence of an integrable function $M$ over $\partial D$ such that
$$\begin{array}{ll}
\displaystyle
\left\vert \frac{\partial}{\partial\nu} v_n(y;x)\right\vert\le M(y),  & \text{a.a. $y\in\partial D$}
\end{array}
$$
since $v_n$ as $n\rightarrow\infty$ highly violates in a neighbourhood of $\sigma$, see \cite{INew}.
So one could not apply the Lebesgue dominated convergence theorem to the last integral on (4.3).

It is open whether the Side B of the singular sources method exits or not
in the present style.
However, in the next subsection we introduce a simple idea 
which enables us to resolve this point.

\subsection{A completely integrated version and the Side B}

In this subsection we reformulate the probe method, that is
introduce another indicator function and show that it is also the indicator function of a reformulated singular sources method.
As a corollary we see this reformulated singular sources method has the Side B.

The idea is to replace the singular solution $G(\,\cdot\,-x)$ with the Green function for the whole domain $\Omega$
and construct the indicator functions of the probe and singular sources methods along the lines indicated by Theorems 1 and 2.

$\quad$

{\bf\noindent Singular solution.}
Given $x\in\Omega$ let $R=R_x\in H^1(\Omega)$ be the weak solution of
$$
\left\{
\begin{array}{ll}
\displaystyle
\Delta R=0,
& 
y\in\Omega,\\
\\
\displaystyle
R=-G(y-x), 
&
\displaystyle
y\in\partial\Omega.
\end{array}
\right.
$$
Here we assume that $\partial\Omega\in C^{1,1}$.  Then we have $R_x\in H^2(\Omega)$.
Define
$$\begin{array}{ll}
\displaystyle
G_{\Omega}(y;x)=G(y-x)+R_x(y), & y\in\Omega\setminus\{x\}.
\end{array}
$$
We see that $G_{\Omega}(\,\cdot\,;x)$ belongs to $H^2$ in a neighbourhood of $\partial\Omega$
and satisfies $G_{\Omega}(\,\cdot\,;x)=0$ on $\partial\Omega$.
The function $G_{\Omega}(\,\cdot\,;x)$ is nothing but the Green function for the domain $\Omega$
and $R_x$ its regular part.

$\quad$

{\bf\noindent  Reflected solution.}
Given $x\in\Omega\setminus\overline{D}$,
let $w=w_x^*(\,\cdot\,)\in H^1(\Omega\setminus\overline{D})$ be the weak solution of (1.4) with
$G(\,\cdot-x)$ replaced with $G_{\Omega}(\,\cdot\,;x)$, that is
$$\left\{\begin{array}{ll}
\displaystyle
\Delta w=0, & y\in\Omega\setminus\overline{D},
\\
\\
\displaystyle
\frac{\partial w}{\partial\nu}=-\frac{\partial}{\partial\nu}G_{\Omega}(y;x), & y\in\partial D,
\\
\\
\displaystyle
w=0, & y\in\partial\Omega.
\end{array}
\right.
\tag {4.4}
$$
Note that we have the expression
$$\begin{array}{ll}
\displaystyle
w_x^*(y)=w_x(y)+z_x(y), 
&
y\in\Omega\setminus\overline{D},
\end{array}
$$
where $w_x(y)$ is the same as Theorem 1 and $z_x(y)=z$ solves
$$\left\{\begin{array}{ll}
\displaystyle
\Delta z=0, & y\in\Omega\setminus\overline{D},
\\
\\
\displaystyle
\frac{\partial z}{\partial\nu}=-\frac{\partial}{\partial\nu}R_x(y), & y\in\partial D,
\\
\\
\displaystyle
z=0, & y\in\partial\Omega.
\end{array}
\right.
$$
Needless to say we have
$$\displaystyle
\forall \epsilon>0\,\sup_{x\in\Omega,\,\text{dist}(x,\partial\Omega)>\epsilon\,}
\,\Vert z_x\Vert_{H^1(\Omega\setminus\overline{D})}<\infty.
$$

$\quad$

The following trivial fact motivates us to introduce another indicator function for the probe method.

\proclaim{\noindent Proposition 4.1.}
Given $\sigma\in N_x$ let $\{v_n(\,\cdot\,;x)\}$ be an arbitrary needle sequence for $(x,\sigma)$,
that is, $v_n(\,\cdot\,;x)$ converges to $G(\,\cdot\,-x)$ in $H^1_{\text{loc}}(\Omega\setminus\sigma)$.
Then the {\it corrected} sequence
$$\displaystyle
\{v_n(\,\cdot\,;x)+R_x\}
$$
converges to $G_{\Omega}(\,\cdot\,;x)$
in $H^1_{\text{loc}}(\Omega\setminus\sigma)$.

\endproclaim

Now, our idea can be explained as follows.
Proposition 4.1 together with formula (1.2) suggests us, in formulating another version of the probe method, instead of the sequence (4.1), to use the sequence
$$\displaystyle
\{<(\Lambda_0-\Lambda_D)((v_n(\,\cdot\,;x)+R_x)\vert_{\partial\Omega}),
(v_n(\,\cdot\,;x)+R_x)\vert_{\partial\Omega}>\}.
$$
Besides, in formulating another version of the singular sources method, considering formula (1.10) together with (1.11), we should use the sequence
$$\displaystyle
\{-<(\Lambda_0-\Lambda_D)((v_n(\,\cdot\,;x)+R_x)\vert_{\partial\Omega}),
(G_{\Omega})_n\vert_{\partial\Omega}>\},
$$
where
$$\begin{array}{ll}\displaystyle
(G_{\Omega})_n(y;x)=G_{\Omega}(y;x)-(v_n(y;x)+R_x(y)), & y\in\Omega\setminus\{x\}.
\end{array}
\tag {4.5}
$$
Then, based on these suggestions, formula (1.3) and  a comparison of the governing equations (1.4) and (4.4) of  $w_x(y)$ and $w_x^*(y)$, respectively,
we should introduce another indicator function for the probe method as follows.

$\quad$

{\bf\noindent Indicator function.}
Define
$$\begin{array}{ll}
\displaystyle
I^*(x)=\int_{\Omega\setminus\overline{D}}\,\vert\nabla w_x^*(y)\vert^2\,dy+\int_D\vert\nabla G_{\Omega}(y;x)\vert^2\,dy,
&
\displaystyle
x\in\Omega\setminus\overline{D}.
\end{array}
$$

$\quad$

Then we obtain the following result.

\proclaim{\noindent  Theorem 7.}

\noindent
(a)   Given $x\in\Omega\setminus\overline{D}$ and $\sigma\in N_x$ let $\{v_n(\,\cdot\,;x)\}$ be an arbitrary needle sequence for $(x,\sigma)$ and $G(\,\cdot\,-x)$.  If $\sigma\cap\overline{D}=\emptyset$, then
we have
$$\displaystyle
\lim_{n\rightarrow\infty}
<(\Lambda_0-\Lambda_D)((v_n(\,\cdot\,;x)+R_x)\vert_{\partial\Omega}),
(v_n(\,\cdot\,;x)+R_x)\vert_{\partial\Omega}>
=I^*(x).
\tag {4.6}
$$

\noindent
(b)  We have, for each $\epsilon>0$
$$\displaystyle
\sup_{x\in\Omega\setminus\overline{D},\,\text{dist}\,(x,\partial D)>\epsilon, 
\text{dist}(x,\partial\Omega)>\epsilon\,}\,I^*(x)<\infty.
$$

\noindent
(c) We have, for each $a\in\partial D$
$$\displaystyle
\lim_{x\rightarrow a}\,I^*(x)=\infty.
$$

\noindent
(d)   For all $x\in\Omega\setminus\overline{D}$ we have 
$$\displaystyle
I^*(x)=w_x^*(x)
\tag {4.7}
$$
and
$$\displaystyle
I^*(x)=I(x)
+2(I^1(x)-w_x^1(x))
+<(\Lambda_0-\Lambda_D)(G(\,\cdot\,-x)\vert_{\partial\Omega}),
G(\,\cdot\,-x)\vert_{\partial\Omega}>.
\tag {4.8}
$$

\endproclaim

{\it\noindent Proof.}
(a) is a consequence of Proposition 4.1 and a similar argument of the derivation of (1.2).

By (1.7) we have
$$\displaystyle
\int_D\vert G_{\Omega}(y;x)\vert^2\,dy\le I(x;\Omega)
\le C\int_D\vert\nabla G_{\Omega}(y;x)\vert^2\,dy,
$$
where $C$ is a positive constant.  Applying the estimate
$$\displaystyle
\forall\epsilon>0\,\sup_{x\in\Omega,\,\text{dist}\,(x,\partial\Omega)>\epsilon}\,\Vert R_x\Vert_{H^1(\Omega)}
<\infty
$$
and the expression $G_{\Omega}(y;x)$, we obtain (b), and (c).

The validity of  (4.7) in (d) is as follows.  
Given $x\in\Omega\setminus\overline{D}$, choose a  $\sigma\in N_x$ with $\sigma\cap\overline{D}=\emptyset$
and $\{v_n(\,\cdot\,;x)\}$ a needle sequence for $(x,\sigma)$ with higher order regularity $v_n(\,\cdot\,;x)\in H^2(\Omega)$.
A similar argument for the derivation of (1.10) we have
$$\displaystyle
-\lim_{n\rightarrow\infty}
<(\Lambda_0-\Lambda_D)((v_n(\,\cdot\,;x)+R_x)\vert_{\partial\Omega}),(G_{\Omega})_n(\,\cdot\,;x)\vert_{\partial\Omega}>
=w_x^*(x),
\tag{4.9}
$$
where the function $(G_{\Omega})_n$ is given  by (4.5).
Here we note that
$$\begin{array}{ll}
\displaystyle
(G_{\Omega})_n(y;x)=G_n(y;x), & y\in\Omega\setminus\{x\}.
\end{array}
$$ 
Besides, since $R_x=-G(\,\cdot\,-x)$ on $\partial\Omega$, we have
$$\begin{array}{ll}
\displaystyle
G_n(y;x)=-(v_n(y;x)+R_x(y)), & y\in\partial\Omega.
\end{array}
\tag {4.10}
$$
Thus this together with (4.9) and (a) yields (4.7).  This is the {\it trick} to connect $I^*(x)$ with $w_x^*(x)$.

Next consider the formula (4.8) in (d).
Again given $x\in\Omega\setminus\overline{D}$, choose a  $\sigma\in N_x$ with $\sigma\cap\overline{D}=\emptyset$
and $\{v_n(\,\cdot\,;x)\}$ a needle sequence for $(x,\sigma)$ with higher order regularity $v_n(\,\cdot\,;x)\in H^2(\Omega)$.
Write
$$\begin{array}{l}
\displaystyle
\,\,\,\,\,\,
<(\Lambda_0-\Lambda_D)((v_n(\,\cdot\,;x)+R_x)\vert_{\partial\Omega}),
(v_n(\,\cdot\,;x)+R_x)\vert_{\partial\Omega}>
\\
\\
\displaystyle
=<(\Lambda_0-\Lambda_D)(v_n(\,\cdot\,;x)\vert_{\partial\Omega}),
(v_n(\,\cdot\,;x)+R_x)\vert_{\partial\Omega}>
+<(\Lambda_0-\Lambda_D)(R_x\vert_{\partial\Omega}),
(v_n(\,\cdot\,;x)+R_x)\vert_{\partial\Omega}>
\\
\\
\displaystyle
=<(\Lambda_0-\Lambda_D)(v_n(\,\cdot\,;x)\vert_{\partial\Omega}),
v_n(\,\cdot\,;x)\vert_{\partial\Omega}>
+<(\Lambda_0-\Lambda_D)(v_n(\,\cdot\,;x)\vert_{\partial\Omega}),
R_x\vert_{\partial\Omega}>
\\
\\
\displaystyle
\,\,\,
+<(\Lambda_0-\Lambda_D)(R_x\vert_{\partial\Omega}),
v_n(\,\cdot\,;x)\vert_{\partial\Omega}>
+<(\Lambda_0-\Lambda_D)(R_x\vert_{\partial\Omega}),
R_x\vert_{\partial\Omega}>
\\
\\
\displaystyle
=<(\Lambda_0-\Lambda_D)(v_n(\,\cdot\,;x)\vert_{\partial\Omega}),
v_n(\,\cdot\,;x)\vert_{\partial\Omega}>
+2<(\Lambda_0-\Lambda_D)(v_n(\,\cdot\,;x)\vert_{\partial\Omega}),
R_x\vert_{\partial\Omega}>
\\
\\
\displaystyle
\,\,\,
+<(\Lambda_0-\Lambda_D)(R_x\vert_{\partial\Omega}),
R_x\vert_{\partial\Omega}>.
\end{array}
$$
Letting $n\rightarrow\infty$, one gets
$$\begin{array}{ll}
\displaystyle
I^*(x) & 
\displaystyle
=I(x) +
2\lim_{n\rightarrow\infty}<(\Lambda_0-\Lambda_D)(v_n(\,\cdot\,;x)\vert_{\partial\Omega}),
R_x\vert_{\partial\Omega}>
\\
\\
\displaystyle
&
\displaystyle
\,\,\,+<(\Lambda_0-\Lambda_D)(R_x\vert_{\partial\Omega}),
R_x\vert_{\partial\Omega}>.
\end{array}
\tag {4.11}
$$
Since $R_x=-G(y-x)$ on $\partial\Omega$, one gets
$$\displaystyle
<(\Lambda_0-\Lambda_D)(R_x\vert_{\partial\Omega}),
R_x\vert_{\partial\Omega}>
=<(\Lambda_0-\Lambda_D)(G(\,\cdot\,-x)\vert_{\partial\Omega}),
G(\,\cdot\,-x)\vert_{\partial\Omega}>.
\tag {4.12}
$$
Besides, we have
$$\displaystyle
<(\Lambda_0-\Lambda_D)(v_n(\,\cdot\,;x)\vert_{\partial\Omega}),
R_x\vert_{\partial\Omega}>
=-<(\Lambda_0-\Lambda_D)(v_n(\,\cdot\,;x)\vert_{\partial\Omega}),
G(\,\cdot\,-x)\vert_{\partial\Omega}>
$$
and by (1.29), one gets
$$
\displaystyle
\lim_{n\rightarrow\infty}<(\Lambda_0-\Lambda_D)(v_n(\,\cdot\,;x)\vert_{\partial\Omega}),
R_x\vert_{\partial\Omega}>
=I^1(x)-w_x^1(x).
$$
Therefore, this together with  (4.11)  and (4.12) yields (4.8).

\noindent
$\Box$

$\quad$

{\bf\noindent Remark 4.1.}
By (1.20) and (1.21) the second term of  (4.8)
is bounded as $x\rightarrow a\in\partial D$
and clearly so is the third term.
Thus we conclude that, as $x\rightarrow\partial D$
$$\displaystyle
I^*(x)\sim_{\partial D} I(x)
$$
and hence
$$\displaystyle
w_x^*(x)\sim_{\partial D} w_x(x).
$$

$\quad$

{\bf\noindent Remark 4.2.}
Using (1.22), we have
$$\displaystyle
I(x)+2(I^1(x)-w_x^1(x))=I(x)+2(I(x)-w_x(x))=3I(x)-2w_x(x).
$$
Applying (1.5) and (1.12) to this right-hand side, from (4.8) one concludes that,
as $x\rightarrow\partial\Omega$
$$
\displaystyle
I^*(x)\sim_{\partial\Omega}\,<(\Lambda_0-\Lambda_D)(G(\,\cdot\,-x)\vert_{\partial\Omega}),G(\,\cdot\,-x)\vert_{\partial\Omega}>.
$$
Note also that, as $x\rightarrow\partial D$
$$\displaystyle
<(\Lambda_0-\Lambda_D)(G(\,\cdot\,-x)\vert_{\partial\Omega}),G(\,\cdot\,-x)\vert_{\partial\Omega}>\sim_{\partial D} 0.
$$
Besides, from (4.7) and (4.8) together with (1.22) we have
$$\displaystyle
w_x^*(x)=w_x(x)+(I^1(x)-w_x^1(x))
+<(\Lambda_0-\Lambda_D)(G(\,\cdot\,-x)\vert_{\partial\Omega}),
G(\,\cdot\,-x)\vert_{\partial\Omega}>.
\tag {4.13}
$$

$\quad$

{\bf\noindent Remark 4.3.}
Note that, for all $x\in\Omega\setminus\overline{D}$, we have
$$\displaystyle
<\Lambda_0(G(\,\cdot\,-x)\vert_{\partial\Omega}), G(\,\cdot\,-x)\vert_{\partial\Omega}>
\not=<\frac{\partial}{\partial\nu}G(\,\cdot\,-x)\vert_{\partial\Omega}, G(\,\cdot\,-x)\vert_{\partial\Omega}>.
$$
The reason is as follows.  Assume that this is not true for a $x\in\Omega\setminus\overline{D}$.  Then, (2.16) yields
$$\displaystyle
<(\Lambda_0-\Lambda_D)(G(\,\cdot\,-x)\vert_{\partial\Omega}), G(\,\cdot\,-x)\vert_{\partial\Omega}><0.
$$
However, this is against the fact
$$\displaystyle
<(\Lambda_0-\Lambda_D)(G(\,\cdot\,-x)\vert_{\partial\Omega}), G(\,\cdot\,-x)\vert_{\partial\Omega}>\ge 0.
$$

$\quad$

{\bf\noindent Remark 4.4.}
By (4.10), the sequence of formula (4.6) of (a) has another expression
$$\begin{array}{l}
\displaystyle
\,\,\,\,\,\,
<(\Lambda_0-\Lambda_D)((v_n(\,\cdot\,;x)+R_x)\vert_{\partial\Omega}),
(v_n(\,\cdot\,;x)+R_x)\vert_{\partial\Omega}>\\
\\
\displaystyle
=<(\Lambda_0-\Lambda_D)(G_n(\,\cdot\,;x)\vert_{\partial\Omega}),
G_n(\,\cdot\,;x)\vert_{\partial\Omega}>.
\end{array}
\tag {4.14}
$$
So, in Theorem 7, we have clarified the limiting behaviour of this {\it third sequence}
provided $\sigma\cap\overline{D}=\emptyset$.  Thus we have obtained the completely integrated version of the Side A of the probe and singular sources methods, where
the common indicator function has two faces $I^*(x)$ and $w_x^*(x)$ and can be calculated by the formula
$$\displaystyle
I^*(x)=w_x^*(x)=\lim_{n\rightarrow\infty}
<(\Lambda_0-\Lambda_D)(G_n(\,\cdot\,;x)\vert_{\partial\Omega}),
G_n(\,\cdot\,;x)\vert_{\partial\Omega}>.
$$

$\quad$

How about the Side B?
By (4.6) and (4.7), we see that the Side B of the singular sources method {\it exists} if one uses the sequence (4.14)
as the indicator sequence, which is nothing but the indicator sequence for the probe method using
the needles $\{v_n(\,\cdot\,;x)+R_x\}$ for $(x,\sigma)$ based on $G_{\Omega}(x;y)$ not $G(\,\cdot\,-x)$.
More pricisely, we have the following result.

\proclaim{\noindent Theorem 8.}  Let $x\in\Omega$ and $\sigma\in N_x$.
Assume that one of two cases (a) and (b) are satisfied:

\noindent
(a) $x\in\overline{D}$;

\noindent
(b) $x\in\Omega\setminus\overline{D}$ and $\sigma\cap D\not=\emptyset$.

\noindent
Then, for any needle sequence $\{v_n(\,\cdot\,;x)\}$ for $(x,\sigma)$
we have 
$$\displaystyle
\lim_{n\rightarrow\infty}
<(\Lambda_0-\Lambda_D)(G_n(\,\cdot\,;x)\vert_{\partial\Omega}),G_n(\,\cdot\,;x)\vert_{\partial\Omega}>=\infty.
$$

\endproclaim

\noindent
Note that the needle sequence added an arbitrary solution in $H^1(\Omega)$ 
of the Laplace equation in $\Omega$, for example $R_x$ in Theorem 7,
has the same blowing up property on the original needle \cite{INew}.  So one can make use of the system
of inequalities  (1.7) described in (ii) of Section 1.2.  Needless to say, like Corollaries 2 and 3 in Section 4, 
we have the Side $B$ and the strong characterization of $\overline{D}$ by using the sequence (4.14).

\proclaim{\noindent Corollary 4.}  A point $x\in\Omega$ belongs to the set $\Omega\setminus\overline{D}$ if and only if there exists a needle $\sigma\in N_x$ and a needle sequence $\{v_n(\,\cdot\,;x)\}$ for $(x,\sigma)$
such that
$$\displaystyle
\sup_n\,<(\Lambda_0-\Lambda_D)(G_n(\,\cdot\,;x)\vert_{\partial\Omega}),
G_n(\,\cdot\,;x)\vert_{\partial\Omega}><\infty.
$$

\endproclaim

$\quad$

And we have an explicit characterization of $\overline{D}$ itself. 

\proclaim{\noindent Corollary 5.}
A point $x\in\Omega$ belongs to the set $\overline{D}$ if and only if
for any needle $\sigma\in N_x$ and needle sequence $\{v_n(\,\cdot\,;x)\}$ for $(x,\sigma)$
we have
$$\displaystyle
\lim_{n\rightarrow\infty}<(\Lambda_0-\Lambda_D)(G_n(\,\cdot\,;x)\vert_{\partial\Omega}),G_n(\,\cdot\,;x)\vert_{\partial\Omega}>=\infty.
$$
\endproclaim

$\quad$

{\bf\noindent Conclusion of this section.}

$\quad$

\noindent
$\bullet$  All the indicator functions $I(x)$, $I^*(x)\equiv w_x^*(x)$, $w_x(x)$ has the same
profile as $x$ approaches a point on $\partial D$.  Their explicit relation is indicated
in the formulae (4.7), (4.8) and (4.13).

\noindent
$\bullet$ The third indicator function like $W_x(x)$ which is different from $I^*(x)=w_x^*(x)$ never appears since 
the Green function vanishes on the outer surface.

\noindent
$\bullet$  The Side B of the singular sources method exists and it is nothing but the Side B of the probe method
using the sequence given by (4.14).  See Corollaries 4 and 5.

\section{Invariance of the blowing up property}

Let $\Omega_1$ and $\Omega_2$ be a pair of bounded domains of $\Bbb R^3$ with smooth boundary
such that $\overline{\Omega}_1\subset\Omega_2$.  Let $D$ be a bounded open set of $\Bbb R^3$
with $C^{1,1}$-boundary
such that $\overline{D}\subset\Omega_1$ and $\Omega_1\setminus\overline{D}$ is connected.
Given $x\in\Bbb R^3\setminus\overline{D}$ let $w=w_x(y;\Omega_j)$ be the solution of
$$\displaystyle
\left\{\begin{array}{ll}
\displaystyle
\Delta w=0, & y\in\Omega_j\setminus\overline{D},
\\
\\
\displaystyle
\frac{\partial w}{\partial\nu}=-\frac{\partial}{\partial\nu}G(y-x), & y\in\partial D,
\\
\\
\displaystyle
w=0, & y\in\partial\Omega_j.
\end{array}
\right.
$$

In this section we show

\proclaim{\noindent Theorem 9.}
We have, for all $x\in\Omega_1\setminus\overline{D}$
$$\displaystyle
w_x(\,\cdot\,;\Omega_2)\vert_{\overline{\Omega_1}\setminus D}-w_x(\,\cdot\,;\Omega_1)
\in C^{0,\frac{1}{2}}(\overline{\Omega_1}\setminus D)
$$
and
$$\displaystyle
\sup_{x\in\Omega_1\setminus\overline{D}}
\Vert 
w_x(\,\cdot\,;\Omega_2)\vert_{\overline{\Omega_1}\setminus D}-w_x(\,\cdot\,;\Omega_1)
\Vert_{C^{0,\frac{1}{2}}
(\overline{\Omega_1}\setminus D)}
<\infty.
$$

\endproclaim

This implies that the blowing up profile of $w_x(x;\Omega_2)$ and $w_x(x;\Omega_1)$ 
as $x\in\Omega_1\setminus\overline{D}$ approaches a point on $\partial D$
coincides with each other.  It means that the blowing up profile
of those indicator functions is invariant
with respect to a monotone perturbation of $\Omega_1$ to $\Omega_2$.
Needles to say, as a conclusion,  the blowing up profile of 
all the indicator functions introduced in this paper is invariant under the perturbation of the outer surface.

\subsection{Proof of Theorem 9}
Let $x\in\Omega_1\setminus\overline{D}$.
Set $z=z_x=w_x(y;\Omega_2)-w_x(y;\Omega_1)$, $y\in\Omega_1\setminus\overline{D}$.
The $z$ satisfies
$$\displaystyle
\left\{\begin{array}{ll}
\displaystyle
\Delta z=0, & y\in\Omega_1\setminus\overline{D},
\\
\\
\displaystyle
\frac{\partial z}{\partial\nu}=0, & y\in\partial D,
\\
\\
\displaystyle
z=w_x(y;\Omega_2), & y\in\partial\Omega_1.
\end{array}
\right.
\tag {5.1}
$$

\proclaim{\noindent Proposition 5.1.}
We have
$$\displaystyle
\sup_{x\in\Omega_2\setminus\overline{D}}\,\Vert z_x\Vert_{H^2(\Omega_1\setminus\overline{D})}<\infty.
\tag {5.2}
$$

\endproclaim
{\it\noindent Proof.}
By Lemma 2.2 in \cite{IR}, we have
$$\displaystyle
\Vert u-v\Vert_{L^2(\Omega_2\setminus\overline{D})}\le C\Vert v\Vert_{L^2(D)},
$$
where $C$ is a positive constant independent of $v$, $v\in H^1(\Omega_2)$ satisfies $\Delta v=0$ in $\Omega_2$
and $u$ solves (1.1) with $\Omega$ and $f$ replaced with $\Omega_2$ and $v$, respectively.
Thus, a limiting procedure yields, for all $x\in\Omega_2\setminus\overline{D}$
$$\displaystyle
\Vert w_x(\,\cdot\,;\Omega_2)\Vert_{L^2(\Omega_2\setminus\overline{D})}
\le C\int_{D}\vert G(y-x)\vert^2\,dy.
$$
Since
$$\displaystyle
\sup_{x\in\Bbb R^3}\int_{D}\frac{dy}{\vert x-y\vert^2}<\infty,
$$
we obtain\footnote{It is possible also to use the estimate (28) on page 214 in \cite{Iwave} for the Helmholtz equation case for sound-soft or sound-hard obstacle case.  However, Lemma 2.2 in \cite{IR} covers more general case.}
$$\displaystyle
\sup_{x\in\Omega_2\setminus\overline{D}}
\Vert w_x(\,\cdot\,;\Omega_2)\Vert_{L^2(\Omega_2\setminus\overline{D})}
\le C'.
\tag {5.3}
$$

Here, by elliptic global regularity \cite{Gr}, we know $w_x(\,\cdot\,;\Omega_2)\in H^2(\Omega_2\setminus\overline{D})$.
Besides, let $U$ be an arbitrary domain of $\Bbb R^3$ such that
$\overline{U}\subset\overline{\Omega_2}$ and $\overline{U}\cap\overline{D}=\emptyset$.
Noting $w_x(\,\cdot\,;\Omega_2)=0$ on $\partial\Omega_2$, one can apply Theorem 9.13 in \cite{GT} to estimate the $H^2(U)$-norm 
of $w_x(\,\cdot\,;\Omega_2)$.  The result is
$$\displaystyle
\Vert w_x(\,\cdot\,;\Omega_2)\Vert_{H^2(U)}
\le C_U\Vert w_x(\,\cdot\,;\Omega_2)\Vert_{L^2(\Omega_2\setminus\overline{D})}.
\tag {5.3}
$$
A combination of  (5.3) and (5.4) yields, for all $x\in\Omega_2\setminus\overline{D}$
$$\displaystyle
\Vert w_x(\,\cdot\,;\Omega_2)\Vert_{H^2(U)}
\le C'_U.
\tag {5.5}
$$
Now choose a special $U$ in such a way that $\partial\Omega_1\subset\overline{U}\subset\overline{\Omega_1}$ and
$\overline{U}\cap\overline{D}=\emptyset$.
Using  the trace theorem and (5.5) with this $U$, we obtain 
$$\displaystyle
\Vert w_x(\,\cdot\,;\Omega_2)\vert_{\partial\Omega_1}\Vert_{H^{\frac{3}{2}}(\partial\Omega_1)}
\le C'\Vert w_x(\,\cdot\,;\Omega_2)\Vert_{H^{2}(U)}\le C'C'_U.
$$
Then the well-posedness of  (5.1) yields (5.2).

\noindent
$\Box$

Thus  a combination of  the Sobolev imbedding and Proposition 5.1 yields Theorem 9.

$\quad$

{\bf\noindent Remark 5.1.}   
Replacing $\Omega_2$ with the general $\Omega$ in the proof of  (5.5),
by virtue of Lemma 2.2 in \cite{IR}
we have, for all $x\in\Omega\setminus\overline{D}$
$$\displaystyle
\Vert w_x(\,\cdot\,;\Omega)\Vert_{H^2(U)}
\le C'_U,
\tag {5.6}
$$
where $U$ satisfies $\overline{U}\subset\overline{\Omega}$ and $\overline{U}\cap\overline{D}=\emptyset$.
By choosing a special $U$ in such a way that $\partial\Omega\subset\overline{U}$ in addition to the property
above and applying the trace theorem together with (5.6), one gets
$$\displaystyle
\sup_{x\in\Omega\setminus\overline{D}}\left\Vert w_x(\,\cdot\,;\Omega)\vert_{\partial\Omega}\right\Vert_{H^{\frac{3}{2}}(\partial\Omega)}
<\infty.
\tag {5.7}
$$
From (5.7) one obtains  (2.18) in Remark 2.1.  The reader should know our approach, that is, introducing
a third indication function is simple and straightforward 
in the sense that we do not make use of nontrivial lemma, that is Lemma 2.2 in \cite{IR}.

A similar comment works also for the article \cite{NPS}.
The author thinks that the approach presented here to obtain (5.7) would work 
in the natural transmission case studied therein with the help of De Giorgi-Nash-Moser theorem,
(e.g.\cite{HL}) as explained in \cite{Itotal}, section 2.2.4. 
Therein an idea of the derivation of  the global boundedness of  the $L^2$-norm of the reflected solution are explained.
However, instead, our third indicator function approach can be also applied.
Those shall be a part of our forthcoming paper.

\section{Final remarks}

Now we have clarified the complete picture of
the probe and singular sources methods applied to a prototype problem
by finding the explicit relationship between
the formulation of the probe and singular sources
methods.  In the picture the needle sequence in the probe method and the notion of
the Carleman function played the role of the bridge 
between both methods.

Constructing the third indicator function, we found that the indicator functions of the both methods
realize two aspects of  this single third indicator function, that is, the energy integral aspect and point-wise one.
A technical advantage of two ways decomposition of this third indicator function  is: the blowing up property of the indicator function for the singular sources
method is almost immediately derived from that of the indicator function for the probe method.
Note that, one can derive the blowing up property of the indicator function for the probe method itself, without making use of an analytical expression of the reflected solution
in a neighbourhood of the surface of the obstacle, see \cite{IProbe}, \cite{Iwave} and \cite{Icrack}  (with Subsection 2.1 in
\cite{Itotal} ) at an early stage and recent article \cite{IR}. 
However, this does not deny the role of  such analytical expression, in some case, it would be useful
for the purpose in clarifying the leading term of the asymptotic profile like (2.7) in Theorem 2.1 of \cite{NPS}.
The point is: we can use both methods to study the blowing up property of the third indicator function, that is
the integration of the both methods.

Now we have shown the {\it skelton} or {\it framework} of the {\it integrated theory} of the probe and singular sources methods.  Needless to say,  if the governing equation or boundary condition are different from those of 
the prototype problem, then we probably
need a modification of the theory or some of statements may not be valid as it is.
However,  the author believes that the third indicator function plays an essential role at the center of the theory.

The one of next directions of this research is to apply this theory to inverse obstacle problems
governed by various partial differential equations including systems, for example, the  Stokes, Maxwell
and Navier systems, and possibly in time domain also.

$$\quad$$

\centerline{{\bf Acknowledgment}}

The author was supported by Grant-in-Aid for
Scientific Research (C)(No. 24K06812) of Japan  Society for
the Promotion of Science.

$$\quad$$

\section{Appendix. A direct proof of (1.18)}

We assume that both $\partial\Omega$ and $\partial D$ are $C^{1,1}$.
Given $x\in\Omega\setminus\overline{D}$ let $w=w_x\in H^2(\Omega\setminus\overline{D})$ be the solution of 
$$\left\{
\begin{array}{ll}
\displaystyle
\Delta w=0, & y\in\Omega\setminus\overline{D},\\
\\
\displaystyle
\frac{\partial w}{\partial\nu}=-\frac{\partial}{\partial\nu}G(y-x), & y\in\partial D,\\
\\
\displaystyle
w=0, & y\in\partial\Omega.
\end{array}
\right.
$$

Let $y\in\Omega\setminus\overline{D}$.
We start with the expression
$$\begin{array}{ll}
\displaystyle
w_x(y)
&
\displaystyle
=\int_{\partial\Omega}
\left(\frac{\partial}{\partial\nu}w_x(z)G(z-y)-
w_x(z)\frac{\partial}{\partial\nu}G(z-y)\right)\,dS(z)\\
\\
\displaystyle
&
\displaystyle
-\int_{\partial D}
\left(\frac{\partial}{\partial\nu}w_x(z)G(z-y)-
w_x(z)\frac{\partial}{\partial\nu}G(z-y)\right)\,dS(z).
\end{array}
\tag {A.1}
$$
Using the boundary conditions on $\partial\Omega$ and $\partial D$, we have
$$\begin{array}{ll}
\displaystyle
w_x(y)
&
\displaystyle
=\int_{\partial\Omega}\frac{\partial}{\partial\nu}w_x(z)G(z-y)\,dS(z)\\
\\
\displaystyle
&
\displaystyle
+\int_{\partial D}
\left(\frac{\partial}{\partial\nu}G(z-x)G(z-y)-
w_x(z)\frac{\partial}{\partial\nu}w_y(z)\right)\,dS(z).
\end{array}
\tag {A.2}
$$
Since $x$ and $y$ outside $D$, we have
$$\begin{array}{ll}
\displaystyle
\int_{\partial D}
\frac{\partial}{\partial\nu}G(z-x)G(z-y)\,dz
&
\displaystyle
=\int_D\Delta G(z-x)G(z-y)\,dz+\int_D\nabla G(z-x)\cdot\nabla G(z-y)\,dz\\
\\
\displaystyle
&
\displaystyle
=\int_D\nabla G(z-x)\cdot\nabla G(z-y)\,dz.
\end{array}
$$
Besides, we have
$$\begin{array}{ll}
\displaystyle
-\int_{\partial D}w_x(z)\frac{\partial}{\partial\nu}w_y(z)\,dS(z)
&
\displaystyle
=\int_{\Omega\setminus\overline{D}}w_x(z)\Delta w_y(z)\,dS(z)
+\int_{\Omega\setminus\overline{D}}\nabla w_x(z)\cdot\nabla w_y(z)\,dz
\\
\\
\displaystyle
&
\displaystyle
=\int_{\Omega\setminus\overline{D}}\nabla w_x(z)\cdot\nabla w_y(z)\,dz.
\end{array}
$$
Thus (A.2) becomes
$$\begin{array}{ll}
\displaystyle
w_x(y)
&
\displaystyle
=\int_{\partial\Omega}\frac{\partial}{\partial\nu}w_x(z)G(z-y)\,dS(z)
\\
\\
\displaystyle
&
\displaystyle
\,\,\,
+\int_{\Omega\setminus\overline{D}}\nabla w_x(z)\cdot\nabla w_y(z)\,dz
+\int_D\nabla G(z-x)\cdot\nabla G(z-y)\,dz.
\end{array}
\tag {A.3}
$$
Next, given $x\in\Omega\setminus\overline{D}$ let $w=w_x^1\in H^2(\Omega\setminus\overline{D})$ be the solution of 
$$\left\{
\begin{array}{ll}
\displaystyle
\Delta w=0, & y\in\Omega\setminus\overline{D},\\
\\
\displaystyle
\frac{\partial w}{\partial\nu}=0, & y\in\partial D,\\
\\
\displaystyle
w=G(y-x), & y\in\partial\Omega.
\end{array}
\right.
$$
Let $y\in\Omega\setminus\overline{D}$.
Same as (A.1) we have the expression
$$\begin{array}{ll}
\displaystyle
w_x^1(y)
&
\displaystyle
=\int_{\partial\Omega}
\left(\frac{\partial}{\partial\nu}w_x^1(z)G(z-y)-
w_x^1(z)\frac{\partial}{\partial\nu}G(z-y)\right)\,dS(z)\\
\\
\displaystyle
&
\displaystyle
\,\,\,
-\int_{\partial D}
\left(\frac{\partial}{\partial\nu}w_x^1(z)G(z-y)-
w_x^1(z)\frac{\partial}{\partial\nu}G(z-y)\right)\,dS(z).
\end{array}
\tag {A.4}
$$
Applying the boundary conditions on $\partial\Omega$ and $\partial D$ to (A.4), we have
$$\begin{array}{ll}
\displaystyle
w_x^1(y)
&
\displaystyle
=\int_{\partial\Omega}
\left(\frac{\partial}{\partial\nu}w_x^1(z)w_y^1(z)-
G(z-x)\frac{\partial}{\partial\nu}G(z-y)\right)\,dS(z)\\
\\
\displaystyle
&
\displaystyle
\,\,\,
+\int_{\partial D}
w_x^1(z)\frac{\partial}{\partial\nu}G(z-y)\,dS(z).
\end{array}
\tag {A.5}
$$
Here we have
$$\begin{array}{ll}
\displaystyle
\int_{\partial\Omega}
\frac{\partial}{\partial\nu}w_x^1(z)w_y^1(z)\,dS(z)
&
\displaystyle
=
\int_{\partial D}
\frac{\partial}{\partial\nu}w_x^1(z)w_y^1(z)\,dS(z)\\
\\
\displaystyle
&
\displaystyle
\,\,\,
+\int_{\Omega\setminus\overline{D}}\Delta w_x^1(z)w_y^1(z)\,dz
+\int_{\Omega\setminus\overline{D}}\nabla w_x^1(z)\cdot\nabla w_y^1(z)\,dz
\\
\\
\displaystyle
&
\displaystyle
=\int_{\Omega\setminus\overline{D}}\nabla w_x^1(z)\cdot\nabla w_y^1(z)\,dz.
\end{array}
$$
Thus (A.5) becomes
$$\begin{array}{ll}
\displaystyle
w_x^1(y)
&
\displaystyle
=\int_{\Omega\setminus\overline{D}}\nabla w_x^1(z)\cdot\nabla w_y^1(z)\,dz
-\int_{\partial\Omega}G(z-x)\frac{\partial}{\partial\nu}G(z-y)\,dS(z)
\\
\\
\displaystyle
&
\displaystyle
\,\,\,
+\int_{\partial D}
w_x^1(z)\frac{\partial}{\partial\nu}G(z-y)\,dS(z).
\end{array}
\tag {A.6}
$$
Here using the boundary condition of $w_y$ on $\partial D$, one has
$$\begin{array}{ll}
\displaystyle
\int_{\partial D}
w_x^1(z)\frac{\partial}{\partial\nu}G(z-y)\,dS(z)
&
\displaystyle
=-\int_{\partial D}
w_x^1(z)\frac{\partial}{\partial\nu}w_y(z)\,dS(z)
\\
\\
\displaystyle
&
\displaystyle
=\int_{\Omega\setminus\overline{D}} w_x^1(z)\Delta w_y(z)\,dz
-\int_{\partial\Omega}w_x^1(z)\frac{\partial}{\partial\nu}w_y(z)\,dS(z)
\\
\\
\displaystyle
&
\displaystyle
\,\,\,
+\int_{\Omega\setminus\overline{D}}\nabla w_x^1(z)\cdot\nabla w_y(z)\,dz\\
\\
\displaystyle
&
\displaystyle
=-\int_{\partial\Omega}G(z-x)\frac{\partial}{\partial\nu}w_y(z)\,dS(z).
\end{array}
$$
Note that we have used
$$\displaystyle
\int_{\Omega\setminus\overline{D}}\nabla w_x^1(z)\cdot\nabla w_y(z)\,dz=0.
$$
Therefore  (A.6) becomes
$$\begin{array}{ll}
\displaystyle
w_x^1(y)
&
\displaystyle
=\int_{\Omega\setminus\overline{D}}\nabla w_x^1(z)\cdot\nabla w_y^1(z)\,dz
-\int_{\partial\Omega}G(z-x)\frac{\partial}{\partial\nu}G(z-y)\,dS(z)
\\
\\
\displaystyle
&
\displaystyle
\,\,\,
-\int_{\partial\Omega}G(z-x)\frac{\partial}{\partial\nu}w_y(z)\,dS(z).
\end{array}
$$
Here we note that
$$\displaystyle
-\int_{\partial\Omega}G(z-x)\frac{\partial}{\partial\nu}G(z-y)\,dS(z)
=\int_{\Bbb R^3\setminus\overline{\Omega}}\nabla G(z-x)\cdot\nabla G(z-y)\,dz.
$$
Thus, we obtain
$$\begin{array}{ll}
\displaystyle
w_x^1(y)
&
\displaystyle
=\int_{\Omega\setminus\overline{D}}\nabla w_x^1(z)\cdot\nabla w_y^1(z)\,dz
+
\int_{\Bbb R^3\setminus\overline{\Omega}}\nabla G(z-x)\cdot\nabla G(z-y)\,dz
\\
\\
\displaystyle
&
\displaystyle
\,\,\,
-\int_{\partial\Omega}G(z-x)\frac{\partial}{\partial\nu}w_y(z)\,dS(z).
\end{array}
\tag {A.7}
$$
A combination of (A.7) and (A.3) together with Definitions 3.1 and 3.2 yields
$$\begin{array}{ll}
\displaystyle
W_x(y)
&
\displaystyle
=I(x,y)+I^1(x,y)
\\
\\
\displaystyle
&
\displaystyle
\,\,\,
+\int_{\partial\Omega}\frac{\partial}{\partial\nu}w_x(z)G(z-y)\,dS(z)
-\int_{\partial\Omega}G(z-x)\frac{\partial}{\partial\nu}w_y(z)\,dS(z).
\end{array}
$$
Therefore, in particular, we obtain
$$\displaystyle
W_x(x)=I(x)+I^1(x).
$$

$\quad$

\end{document}